\documentclass[sn-mathphys]{sn-jnl}

\usepackage[utf8]{inputenc}
\usepackage[english]{babel}
\usepackage{amsthm}
\usepackage{amsmath}
\usepackage{amssymb}
\usepackage{amsfonts}
\usepackage{dsfont}
\usepackage{mathrsfs}
\usepackage{graphicx}

\usepackage{enumerate}
\usepackage{hyperref}
\usepackage{indentfirst}
\usepackage{color}

\newcommand{\dist}[2]{\langle {#1} , {#2} \rangle}

\newcommand{\norm}[2]{\vert  {#1}  \vert_{#2}}
\newcommand{\norma}[2]{\|  {#1} \|_{#2}}

\newcommand{\D}[1]{\displaystyle{#1}}

\usepackage{verbatim}

\newtheorem{theorem}{Theorem}[section]
\newtheorem{lemma}[theorem]{Lemma}

\theoremstyle{remark}
\newtheorem{definition}[theorem]{\bf Definition}
\newtheorem{remark}[theorem]{\bf Remark}
\newtheorem{hypothesis}{\bf Hypothesis}

\begin{document}

\title[Analysis and Optimal Control of Chemotaxis-Consumption Models]{A Review on the Analysis and Optimal Control of Chemotaxis-Consumption Models}


\author*[1]{\fnm{André Luiz} \sur{Corrêa Vianna Filho}}\email{viannafilhoandre@gmail.com}
\equalcont{These authors contributed equally to this work.}

\author[1]{\fnm{Francisco} \sur{Guillén-González}}\email{guillen@us.es}
\equalcont{These authors contributed equally to this work.}

\affil[1]{\orgdiv{EDAN}, \orgname{Universidad de Sevilla}, \orgaddress{\street{Av. Reina Mercedes}, \city{Sevilla}, \postcode{41013}, \state{Sevilla}, \country{Spain}}}


\abstract{In the present review we focus on the chemotaxis-consumption model $\partial_t u - \Delta u  = -  \nabla \cdot (u \nabla v)$ and $\partial_t v - \Delta v  = - u^s v$ in $(0,T) \times \Omega$, for any fixed $s \geq 1$, endowed with isolated boundary conditions and nonnegative initial conditions, where $(u,v)$ model cell density and chemical signal concentration. Our objective is to present an overview of the related literature and latest results on the aforementioned model concerning the following three distinct research lines we have obtained in \cite{ViannaGuillen2023uniform,guillen2023convergence,guillen2023optimal,guillen2022optimal}: the mathematical analysis, the numerical analysis and the related optimal control theory with a bilinear control acting on the chemical equation.}




\maketitle


\section{Introduction}\label{sec1}

In microbiology, chemotaxis is understood as the directed migration of cells in response to a concentration gradient of a certain chemical substance, either toward attractant chemicals or away from repelents \cite{STOCK200971}. Chemotaxis plays an essential role in many biological processes such as wound healing, the immune cells migration, the migration of bacteria, among others. It is also an important factor in some undesired events such as tumor growth, cancer metastasis and inflamatory diseases, see for example \cite{murphy2001chemokines,wang2011signaling}. In unicellular organisms, chemotaxis is frequently related to the search for nutrients \cite{FRANCAKOH20101705} and there are studies on its applications to the degradation of polluting substances \cite{pandey2002bacterial,parales2000toluene}.

  The introduction of one of the first mathematical models for chemotaxis is attributed Keller and Segel in two works from 1970 and 1971 \cite{keller1970initiation,keller1971model} which are also regarded by some authors as a development of the work of Patlak \cite{patlak1953random}. Since then, the research on this topic gave rise to different related models such as models with chemoattraction or chemorepulsion, combined with production or consumption of the chemical substance, with the presence of a logistic growth of the population of cells, models for angiogenesis, haptotaxis and so on, covering a wide variety of applications of practical interest. From the mathematical point of view, the aforementioned models possess interesting and challenging features that attracted the attention of many authors along the years and make these models still relevant nowadays \cite{bellomo2015toward,horstmann20031970,horstmann20041970}.

  In the present review, we focus on a model that describes a situation where, inside an isolated, bounded and connected region $\Omega$ of the $d$-dimensional space $\mathbf{R}^d$ ($d = 1, 2, 3$) with boundary $\Gamma$, the cells diffuse and are attracted by the chemical substance that, in its turn, diffuses and is consumed by the cells. Let $u=u(t,x)$ and $v=v(t,x)$ be the density of cell population and the concentration of chemical substance, respectively, on each spatial point  $x \in \Omega$ and time $t \geq 0$. This model is governed by the initial-boundary PDE problem
  \begin{equation}
    \left\{\begin{array}{l}
      \partial_t u - \Delta u  = - \nabla \cdot (u \nabla v), \quad
      \partial_t v - \Delta v  = - u^s v, 
      \quad \hbox{in $(0,+\infty)\times \Omega$,}
\\ [6pt]
      \partial_{\bf n} u \vert_{\Gamma}  =  \partial_{\bf n} v \vert_{\Gamma} = 0, \quad
      u(0)  = u^0, \quad v(0) = v^0,
      \quad \hbox{in $\Omega$,}
    \end{array}\right.
    \label{problema_P} 
  \end{equation}
  where $\nabla \cdot (u \nabla v)$ is the chemotaxis term, $u^s v$ is the consumption term, with $s \geq 1$, and the normal derivative of $u$ on the boundary is denoted by $\partial_{\bf n} u$. We assume that the initial conditions $u^0 = u^0(x)$ and $v^0 = v^0(x)$ are nonnegative given functions.

  We aim for presenting an overview of the latest results on the chemotaxis consumption model \eqref{problema_P} concerning the following three distinct research lines: (i) the mathematical analysis, (ii) the numerical analysis and (iii) the related optimal control theory. The controlled model consists of the chemotaxis-consumption model \eqref{problema_P} with a control $f: (0,T) \times \Omega \rightarrow \mathbb{R}$, being $T > 0$ a fixed and finite final time, acting directly on the chemical equation through the bilinear term $f v 1_{\Omega_c}$:
  \begin{equation} \label{problema_P_controlado}
    \left\{\begin{array}{l}
      \partial_t u - \Delta u  = - \nabla \cdot (u \nabla v), \quad
      \partial_t v - \Delta v  = - u^s v + f v 1_{\Omega_c},
      \  \hbox{in $(0,T)\times \Omega$,} \\ [6pt]
      \partial_{\bf n} u  \vert_{\Gamma}  =  \partial_{\bf n} v  \vert_{\Gamma} = 0, \quad
      u(0)  = u^0, \quad v(0) = v^0,
      \  \hbox{in $\Omega$,}
    \end{array}\right.
  \end{equation}
  where $\Omega_c \subset \Omega$ is the control domain and $1_{\Omega_c}$ is its characteristic function. The control of chemotaxis systems through the distributed action on the chemical equation has been considered in previous studies \cite{ryu2001optimal,guillen2020optimal,guillen2020bi,guillen2020regularity,lopez2021optimal,silva2022bilinear,tang2022optimal,guillen2023optimal}. In contrast to \cite{ryu2001optimal}, where the control acts as a linear term, the advantage of using a bilinear term such as $fv1_{\Omega_c}$ to control the system is that it allows us to preserve the positivity of $v$ independently of the sign of $f$.

    In what follows, we begin with an overview on the mathematical analysis of the chemotaxis-consumption models \eqref{problema_P}. Then we turn to the numerical analysis and finally to the optimal control theory of these models.


  \section{Analytical Results} \label{sec:Analytical}

    First we recall some developments in the chemotaxis model \eqref{problema_P} for the most studied case, with $s = 1$ (hence the consumption term $-uv$ is bilinear). In this case, the available results are reached through the use of classical-in-time solution tools (Amman's theory \cite{amann1995linear}), considering PDEs with constant coefficients and smooth domains. After that, we turn to the latest results considering $s \geq 1$, more general domains and a setting which is more appropriate to approach the time approximation and the optimal control of \eqref{problema_P_controlado}. 
    
    In \cite{tao2012eventual}, existence of global weak solutions which become smooth after a sufficiently large period of time is proved in smooth and convex $3$D domains. More recently, a parabolic-elliptic simplification of \eqref{problema_P}, also for $s = 1$, is studied in \cite{tao2019global}, yielding results on the existence and long-time behavior of global classical solutions in $d$-dimensional smooth domains.
    
    Still considering $s = 1$, there are some studies on the coupling of \eqref{problema_P} with models for incompressible fluids. In \cite{lorz2010coupled}, the author proves local existence of weak solutions for the chemotaxis-Navier-Stokes equations in $3$D smooth domains, while in \cite{duan2010global} the existence of global classical solutions is attained near constant states. In \cite{winkler2012global}, considering smooth and convex domains, existence and uniqueness of a global classical solution for the chemotaxis-Navier-Stokes equations is proved in $2$D and existence of global weak solutions which become smooth after a large enough period of time is proved for the chemotaxis-Stokes equations in $3$D. In \cite{jiang2015global} the results of \cite{winkler2012global} on the existence of solution are extended to non-convex domains, but we remark that some estimates that were time-independent in \cite{winkler2012global} become time-dependent in \cite{jiang2015global}. In \cite{winkler2014stabilization} the author studies the assymptotic behavior of the chemotaxis-Navier-Stokes equations in 2D domains with the chemotaxis and consumption terms generalized by using adequate functions depending on the chemical substance, proving the convergence towards constant states in the $L^{\infty}$-norm. Finally, in \cite{winkler2016global} existence of global weak solutions for the chemotaxis-Navier-Stokes equations is established in $3$D smooth and convex domains and in \cite{winkler2017far} the assymptotic behavior of these solutions is studied.

    An interesting and challenging feature of chemotaxis models, both from theoretical and numerical point of view, is that the $L^{\infty}$-norm of the cell density $u$ may blow up in finite time. Some studies focus on the proof of the existence of blowing-up solutions, see for instance \cite[Theorem 3.3]{bellomo2015toward}, while others are dedicated to the proof of existence of uniformly in time  bounded solutions. When it comes to the model \eqref{problema_P}, with $s = 1$, this question has been addressed for $2$D smooth and convex domains, because existence and uniqueness of classical and uniformly bounded solutions is proved in \cite{tao2012eventual}. On the other hand, as far as we know, this question remains open for $3$D domains.
    
    Studying conditions that could lead to no-blow-up results for \eqref{problema_P}, with $s = 1$, some researchers advanced under the assumption of adequate constraints on $\norma{v^0}{L^{\infty}(\Omega)}$ with respect to chemotaxis coefficient. On this subject, we refer the interested reader to \cite{baghaei2017boundedness} and \cite{tao2011boundedness}, for the problem \eqref{problema_P} with $s = 1$. In addition, we also have \cite{frassu2021boundedness} and \cite{fuest2019analysis}, where these results are extended to other related chemotaxis models with consumption.

    For an exhaustive review on the analytical results on the model \eqref{problema_P} (for $s = 1$) and some variants we refer the reader to the recent survey \cite{lankeit2023depleting}, which includes great part of the studies cited above.

    We remark that the aforementioned works were carried out using classical-in-time solution tools, through the heat semigroup theory. The Stokes semigroup is also used in the case of models coupled with fluids. The theoretical background that is applied relies on the fact that the diffusion operator is the Laplacian operator and requires smooth domains and smooth coefficients.
    
    Taking that into account, the previous classical-in-time theory is not well suited neither to the numerical approximation of \eqref{problema_P} nor to the study of optimal control problems subject to the controlled problem \eqref{problema_P_controlado}. Indeed, with respect to the numerical approximation, one usually employs a weak formulation of the problem posed in more general domains. The controlled problem \eqref{problema_P_controlado}, in its turn, contains the control term $f v$ where $f = f(t,x)$ is usually only a $L^q$-function, which can be seen as a non-smooth coefficient.

    The facts that have been exposed so far motivated the extension of the classical theory about \eqref{problema_P} to the weak framework in \cite{ViannaGuillen2023uniform}, which is more suitable to the design of stable and convergent (time) discrete schemes and to the study of optimal control problems subject to \eqref{problema_P_controlado}. This lead to results about the existence and regularity of solutions of \eqref{problema_P} in a weak setting for a larger class of considered domains, avoiding the convexity or smoothness assumptions on $\Omega$.
    
    In addition, a more general consumption term given by $- u^s v$, for $s \geq 1$, is considered in \cite{ViannaGuillen2023uniform}. Actually, we could understand the consumption term as being of the form $-g(u) v$, where we may take $g(u)$ as being any sufficiently regular function generalizing the potential behavior of the prototype $g(u) = u^s$, for $s \geq 1$. This kind of generalization may be interesting from the modeling point of view. Indeed, different assumptions during the modeling process can lead to different terms in the equations and, therefore, it would be useful to have results for  models with more general terms, rather than for a specific case.
    
    It is also interesting to know, what is the effect of distinct chemotaxis and consumption terms with respect to the properties of the solutions, such as existence, uniqueness, regularity, boundedness, asymptotic behavior and so on. For instance, in \cite{ViannaGuillen2023uniform}, the effects of the consumption power $s$ in the regularity of the solutions has  been addressed.

    We would also like to make a comment regarding the rigor of the calculations. We have observed that, in some papers on analysis of chemotaxis models, singular functions are taken as test function, as $log(u)$ for instance. In the context of works such as \cite{winkler2012global} and \cite{tao2012eventual}, where, first of all, the authors prove the existence of local-in-time classical solutions, one could verify if it is possible to apply any strong ``maximum principle'' result to prove that $u$ is strictly positive in the whole domain and then functions of the cell density $u$ with singularity at zero can be used as test functions. Otherwise, according to our understanding, one should take actions to guarantee that all the computations carried out are rigorous. To illustrate it, we refer the reader to the techniques used in \cite{ViannaGuillen2023uniform} to the chemotaxis-consumption models \eqref{problema_P} and in \cite{jungel2022analysis}, where a cross-diffusion model is studied and, to make rigorous computations, some regularization procedures are used, such as taking $log(u+\epsilon)$ as a test function instead of $log(u)$.
    
    Since the review \cite{lankeit2023depleting} accounts for most of the studies related to \eqref{problema_P} with $s = 1$ cited above, we proceed and finish this section exploring more the contributions and main results of \cite{ViannaGuillen2023uniform}, where $s \geq 1$ is considered. To carry out with their analysis, the authors introduced a regularization process by using truncated models depending on a truncation parameter $m \in \mathbb{N}$:
    \begin{equation}
    \left\{
      \begin{array}{l}
        \partial_t u_m - \Delta u_m  = - \nabla \cdot (T^m(u_m) \nabla v_m), \quad
        \partial_t v_m - \Delta v_m  = - T^m(u_m)^s v_m, \\ [6pt]
        \partial_{\bf n} u_m \vert_{\Gamma}  = \partial_{\bf n} v_m \vert_{\Gamma} = 0, \quad
        u_m(0)  = u^0_m, \quad v_m(0) = v^0_m,
      \end{array}
      \right.
      \label{problema_P_m_intro}
    \end{equation}
    where $u^0_m \geq 0$ and $v^0_m \geq 0$ are suitable regular approximations of $u^0$ and $v^0$, respectively, and $T^m(\cdot)$ is a suitable truncation of the identity function (bounded from above and from below) defined by
    \begin{equation} \label{truncamento_limitado_da_identidade}
      T^m(u) = \left \{
      \begin{array}{cl}
        - 1, & \mbox{ if } u \leq -2,  \\
        C^2 \mbox{ extension}, & \mbox{ if } u \in (-2,0), \\
        u, & \mbox{ if } u \in [0, m], \\
        C^2 \mbox{ extension}, & \mbox{ if } u \in (m,m + 2), \\
        m + 1, & \mbox{ if } u \geq m + 2.
      \end{array}
      \right .
    \end{equation}
    
    These truncated models \eqref{problema_P_m_intro} are easier to analyze both from the theoretical and numerical points of view and it is proved in \cite{ViannaGuillen2023uniform} that the solutions of the truncated models converge to weak solutions of \eqref{problema_P} as $m \to \infty$.
        
    With the objective of enlarging the class of considered domains, the results are stated and demonstrated in terms of the regularity of the Poisson-Neumann problem
    \begin{equation} \label{Neumann_problem} 
      \left \{ \begin{array}{rl}
        - \Delta z + z & = h  \, \mbox{ in } \Omega, \\
        \partial_{\bf n} z \vert_{\Gamma} & = 0 \, \mbox{ on } \Gamma,
      \end{array} \right.
    \end{equation}
    and, when necessary, in terms of a technical hypothesis, both stated below.

\

    \begin{definition}[\bf Regularity of the Poisson-Neumann problem] \label{defi_regularidade_H_m}
    Let $z \in H^1(\Omega)$ is a weak solution of \eqref{Neumann_problem} with $h \in L^p(\Omega)$. If this implies that $z \in W^{2,p}(\Omega)$ with
    \begin{equation*}
      \norma{z}{W^{2,p}(\Omega)} \leq C \norma{h}{L^p(\Omega)}, 
    \end{equation*}  
    then we say that the Poisson-Neumann problem \eqref{Neumann_problem} has the $W^{2,p}$-regularity. In the hilbertian case $p = 2$ we say $H^2$-regularity.
    \hfill $\square$
  \end{definition}

\
  
  \begin{remark}   
    According to Grisvard \cite{Grisvard}, if $h \in L^p(\Omega)$, $p \in [1,\infty]$, and the boundary $\Gamma$ is at least $C^{1,1}$, then the Neumann problem \eqref{Neumann_problem} has the $W^{2,p}$-regularity for all $p \in [1,\infty]$. This result is also true if $\Omega$ is a polygon, that is, a polyhedron in $\mathbb{R}^2$, or if $\Omega$ is convex and $p=2$.
    \hfill $\square$
  \end{remark}

\

    \begin{hypothesis}
      For each $z \in H^2(\Omega)$ such that $\partial_{\eta} z \vert_{\Gamma} = 0$ there is a sequence $\{ \rho_n \} \subset C^2(\overline{\Omega})$ such that $\partial_{\eta} \rho_n \vert_{\Gamma} = 0$ and $\rho_n \to z$ in $H^2(\Omega)$.
      \label{hypothesis_density}
    \end{hypothesis}

\

    \begin{remark}
      In order to show that the Hypothesis \ref{hypothesis_density} is not too restrictive, we prove in \cite[Lemma 34]{ViannaGuillen2023uniform} that Hypothesis \ref{hypothesis_density} is satisfied if the Poisson-Neumann problem has the $W^{3,p}$-regularity (see definition \ref{defi_regularidade_H_m} above), for $p > d$. This is true, in particular, if $\Gamma$ is $C^{2,1}$ (see \cite{Grisvard}).
      \hfill $\square$
    \end{remark}    
   
   \

    We are in position to highlight the main results of \cite{ViannaGuillen2023uniform}.
       
    \begin{theorem}[\bf $\boldsymbol{3}$D. Existence of global weak solutions of \eqref{problema_P}] \label{theorem_3D_existencia}
      Let $\Omega \subset \mathbb{R}^3$ be a bounded domain such that the Neumann problem \eqref{Neumann_problem} has the $H^2$-regularity (see definition \ref{defi_regularidade_H_m}) and Hypothesis \ref{hypothesis_density} is satisfied. Let $u^0 \in L^{1 + \varepsilon}(\Omega)$, for some $\varepsilon > 0$, if $s = 1$, and $u^0 \in L^s(\Omega)$, if $s > 1$, and $v^0 \in H^1(\Omega) \cap L^{\infty}(\Omega)$ be non-negative functions. Then there is a non-negative weak solution $(u,v)$ of the original problem \eqref{problema_P}, for $s \geq 1$, obtained through a limit of non-negative solutions $(u_m,v_m)$ of the regularized problems \eqref{problema_P_m_intro} as $m \to \infty$ and such that
      \begin{equation}
        \left \{
        \begin{array}{rl}
          \D{\int_{\Omega}}{u(t,x) \ dx} = \D{\int_{\Omega}}{u^0(x) \ dx}, \ a.e. \ t \in (0,\infty) \\
          0 \leq v(t,x) \leq \norma{v^0}{L^{\infty}(\Omega)}, \ a.e. \ (t,x) \in (0,\infty) \times \Omega,
        \end{array}
        \right.
        \label{pointwise_properties_3D}
      \end{equation}
      \begin{equation*}
        \left \{
        \begin{array}{rl}
          u \in L^{\infty}(0,\infty;L^s(\Omega)) \cap L^{5s/3}_{loc}([0,\infty);L^{5s/3}(\Omega)), \mbox{ if } s \geq 1, \\
          u^{s/2} \nabla v \in L^2(0,\infty;L^2(\Omega)), \mbox{ if } s \geq 1,
        \end{array}
        \right.
      \end{equation*}
      \begin{equation*}
        \left \{
        \begin{array}{rl}
          \nabla u \in L^2(0,\infty;L^s(\Omega)) \cap L^{\frac{5s}{3+s}}_{loc}([0,\infty);L^{\frac{5s}{3+s}}(\Omega)), & \mbox{ if } s \in [1,2), \\
          \nabla u \in L^2(0,\infty;L^2(\Omega)), & \mbox{ if } s \geq 2,
        \end{array}
        \right.
      \end{equation*}
      \begin{equation*}
        \left \{
        \begin{array}{rl}
          u \nabla v \in L^2(0,\infty;L^s(\Omega)), & \mbox{ if } s \in [1,2), \\
          u \nabla v \in L^2(0,\infty;L^2(\Omega)), & \mbox{ if } s \geq 2
        \end{array}
        \right.
      \end{equation*}
      and
      \begin{equation*}
        v \in L^{\infty}(0,\infty;H^1(\Omega)) \cap L^2(0,\infty;H^2(\Omega)),
        \quad
        \nabla v \in L^4(0,\infty; L^4(\Omega)).
      \end{equation*}
    \end{theorem}
    
    \begin{remark} \label{remark_initial_conditions_weak_solutions}
      We remark that, from the regularities of $u$ and $v$ that are listed in Theorem \ref{theorem_3D_existencia}, we can conclude that
      \begin{equation*}
        \left \{
        \begin{array}{rl}
          u_t \in L^2  \big(0,\infty;  W^{1,s/(s-1)}(\Omega) '  \big), & \mbox{ if } s \in [1,2), \\
          u_t \in L^2  \big( 0,\infty; H^1(\Omega) '  \big), & \mbox{ if } s \geq 2,
        \end{array}
        \right.
      \end{equation*}
      where $ W^{1,s/(s-1)}(\Omega) ' $ denotes the dual space of  $W^{1,s/(s-1)}(\Omega)$ (and the same for $H^1(\Omega) ' $)
      and
      \begin{equation*}
        v_t \in L^2(0,\infty;L^{3/2}(\Omega)).
      \end{equation*}
      Attending to the regularity of $(u,v)$ given so far, one has that $(u,v)$ satisfies the $u$-equation of \eqref{problema_P} in a variational sense and the $v$-equation pointwisely $a.e.$ in $(0,\infty) \times \Omega$. Moreover, the initial conditions have a sense because, thanks to the regularity of $u$, $v$, $u_t$ and $v_t$ given above, one has that $(u,v)$ is weakly continuous from $[0,\infty)$ to $L^s(\Omega) \times H^1(\Omega)$, if $s \in [1,2]$, and $L^2(\Omega) \times H^1(\Omega)$, if $s \geq 2$ (owing to a result given, for instance, in Chapter 3 of \cite{Temam}).
      \hfill $\square$
    \end{remark}

\
    
    \begin{remark}
      Note that, for $s \in [1,2]$, the regularity of the fluxes of the $u$-equation of \eqref{problema_P}, namely, self diffusion $\nabla u$ and chemotaxis $u \nabla v$, increase as $s$ increases. But, when we consider $s > 2$, the regularity of $\nabla u$ and $u \nabla v$ does not increase as $s$ increases anymore. On the other hand, the regularity of the chemical variable $v$ is independent of $s$.
      \hfill $\square$
    \end{remark}
    
    \begin{theorem}[\bf $2$D. Existence and uniqueness of global strong solution] \label{theorem_2D_existencia_unicidade_estimativas}
      Let $\Omega \subset \mathbb{R}^2$ be a bounded domain such that the Neumann problem \eqref{Neumann_problem} has the $W^{2,3}$-regularity (see definition \ref{defi_regularidade_H_m}) and Hypothesis \ref{hypothesis_density} is satisfied. Let $u^0 \in H^2(\Omega)$ and $v^0 \in H^2(\Omega)$ be such that $u^0 \geq 0$ and $v^0 \geq 0$ in $\Omega$. Then there is a unique non-negative solution $(u,v)$ for the original problem \eqref{problema_P}, for $s \geq 1$, satisfying  \eqref{pointwise_properties_3D} and the regularity
      \begin{equation*}
        \begin{array}{c}
          u, v \in L^{\infty}(0,\infty;H^2(\Omega)), \quad 
          \Delta u, \Delta v, u_t, v_t \in L^2(0,\infty;H^1(\Omega)).
        \end{array}
      \end{equation*}
      In particular, $u$ does not blow-up neither at finite nor infinite time, that is, $u \in L^{\infty}(0,\infty;L^{\infty}(\Omega))$ (recall that in Theorem \ref{theorem_3D_existencia} we already have $v \in L^{\infty}(0,\infty;L^{\infty}(\Omega))$). Consequently, there is $m_0 \in \mathbb{N}$ such that, for all $m \in [m_0, \infty)$, the solution of \eqref{problema_P_m_intro} is also the solution of \eqref{problema_P}, that is,
      \begin{equation*}
        (u_m,v_m) = (u,v) \quad a.e. \mbox{ in } (0,\infty) \times \Omega.
      \end{equation*}
    \end{theorem}
    
    In this case, both equations of \eqref{problema_P} are satisfied $a.e.$ in $(t,x) \in (0,\infty) \times \Omega$.
    
    \

    In the rest of this section, we give an idea of how Theorems \ref{theorem_3D_existencia} and \ref{theorem_2D_existencia_unicidade_estimativas} are proved in \cite{ViannaGuillen2023uniform}. To show the existence of solutions of \eqref{problema_P} as limits of solutions of \eqref{problema_P_m_intro}, as it is stated in Theorems \ref{theorem_3D_existencia} and \ref{theorem_2D_existencia_unicidade_estimativas}, it is necessary to pass to the limit in \eqref{problema_P_m_intro} as $m \to \infty$. To this end, we need some $m$-independent estimates for $(u_m,v_m)$.

    In \cite{ViannaGuillen2023uniform}, when the truncated model \eqref{problema_P_m_intro} is studied, some direct $m$-independent estimates for $(u_m,v_m)$ are proved. We have $m$-independent estimates for $u_m$ in $L^{\infty}(0,\infty;L^1(\Omega))$ and for $v_m$ in $L^{\infty}(0,\infty;L^{\infty}(\Omega))$ and for $\nabla v_m$ in $L^{\infty}(0,\infty;L^2(\Omega))$. However, these estimates are not sufficient to pass to the limit in \eqref{problema_P_m_intro} as $m \to \infty$ and stronger estimates are necessary.

    The basic idea to obtain additional \emph{a priori} $m$-independent estimates is that the effects of the consumption and chemotaxis terms cancel. It is useful to consider some formal calculations to illustrate how it works. Suppose $(u,v)$ is a regular enough solution to the original problem \eqref{problema_P} with $u,v > 0$. Consider the change of variable $z = \sqrt{v}$, then  problem \eqref{problema_P} can be rewritten as
    \begin{equation}
      \begin{array}{rl}
        \partial_t u - \Delta u & = - \nabla \cdot (u \nabla (z)^2 ), \\
        \partial_t z - \Delta z - \dfrac{\norm{\nabla z}{}^2}{z} & = - \dfrac{u^s z}{2} ,\\
        \partial_{\bf n} u \vert_{\Gamma} & =  \partial_{\bf n} z \vert_{\Gamma} = 0, \\
        u(0) & = u^0, \quad z(0) = \sqrt{v^0}.
      \end{array}
      \label{problema_P_u_z} \\
    \end{equation}
    We are going to obtain estimates for $u$ and $z$ and then extract estimates for $v$ from the estimates of $z$. For this, we consider a function $g(u)$ such that $g''(u) = u^{s-2}$. Formally, assuming $u,z > 0$ we can use 
    \begin{equation}
      g'(u) = \left \{
      \begin{array}{rl}
        \dfrac{u^{s-1}}{(s-1)} &, \mbox{ if } s > 1, \\[6pt]
        ln(u) &, \mbox{ if } s = 1.
      \end{array}
      \right.
      \label{funcao_teste_formal}
    \end{equation}
    as a test function in the $u$-equation of \eqref{problema_P_u_z}, obtaining
    \begin{equation*}
      \dfrac{d}{dt} \int_{\Omega}{g(u) \ dx} + \int_{\Omega}{g''(u) \norm{\nabla u}{}^2 \ dx} = \int_{\Omega}{u g''(u) \nabla (z^2) \cdot \nabla u}
    \end{equation*}
    and, since $u g''(u) = u^{s-1}$, we have
    \begin{equation} \label{estimativa_formal_u}
      \begin{array}{rl}
        \dfrac{d}{dt} \D{\int_{\Omega}{g(u) \ dx} + \int_{\Omega}{u^{s-2} \norm{\nabla u}{}^2 \ dx}} & = \D{\int_{\Omega}{u^{s-1} \nabla (z^2) \cdot \nabla u \ dx}} \\
        & = \dfrac{1}{s} \D{\int_{\Omega}{\nabla (z^2) \cdot \nabla (u^s) \ dx}}.
      \end{array}
    \end{equation}
    On the other hand, we can test the $z$-equation of \eqref{problema_P_u_z} by $- \Delta z$. Then we obtain
    \begin{equation} \label{estimativa_formal_z}
      \begin{array}{rl}
        \dfrac{1}{2} \dfrac{d}{dt} \D{\int_{\Omega}}{\norm{\nabla z}{}^2 \ dx} & + \D{\int_{\Omega}}{\norm{\Delta z}{}^2 \ dx} + \D{\int_{\Omega}}{\dfrac{\norm{\nabla z}{}^2}{z} \Delta z \ dx} + \dfrac{1}{2} \D{\int_{\Omega}}{u^s \norm{\nabla z}{}^2 \ dx} \\[12pt] 
        & = - \dfrac{1}{4} \D{\int_{\Omega}}{\nabla (u^s) \cdot \nabla (z^2) \ dx}.
      \end{array}
    \end{equation}
    Hence, if we add \eqref{estimativa_formal_z} to $s/4$ times \eqref{estimativa_formal_u}, then the two terms on the right hand side cancel each other and we obtain the time differential equation
    \begin{equation} \label{eq_formal_u_w}
      \begin{array}{rl}
      &  \dfrac{d}{dt} \left[ \dfrac{s}{4} \D{\int_{\Omega}{g(u) \ dx}}  + \dfrac{1}{2} \D{\int_{\Omega}{\norm{\nabla z}{}^2 \ dx}} \right] + \dfrac{s}{4} \D{\int_{\Omega}{u^{s-2} \norm{\nabla u}{}^2 \ dx}}  \\[12pt]
        & \qquad + \dfrac{1}{2} \D{\int_{\Omega}{u^s \norm{\nabla z}{}^2 \ dx} + \D{\int_{\Omega}{\norm{\Delta z}{}^2 \ dx}} + \int_{\Omega}{\dfrac{\norm{\nabla z}{}^2}{z} \Delta z \ dx}} = 0.
      \end{array}
    \end{equation}
    The main idea now is to estimate, from below, the term
    \begin{equation} \label{termo_problematico}
      \D{\int_{\Omega}{\norm{\Delta z}{}^2 \ dx}} + \int_{\Omega}{\dfrac{\norm{\nabla z}{}^2}{z} \Delta z \ dx},
    \end{equation}
    for which the following lemma is essential.
    
    \begin{lemma}[\bf \cite{ViannaGuillen2023uniform}] \label{lemma_termo_fonte_final}
      Suppose that the Poisson-Neumann problem \eqref{Neumann_problem} has the $H^2$-regularity and assume that Hypothesis \ref{hypothesis_density} holds. Then there exist positive constants $C_1, C_2 > 0$ such that
      \begin{align*}
        & \int_{\Omega}{\norm{\Delta z}{}^2 \ dx} + \int_{\Omega}{\frac{\norm{\nabla z}{}^2}{z} \Delta z \ dx} \geq \\
        & C_1 \Big ( \int_{\Omega}{\norm{D^2 z}{}^2 \ dx} + \int_{\Omega}{\frac{\norm{\nabla z}{}^4}{z^2} \ dx} \Big ) - C_2 \int_{\Omega}{\norm{\nabla z}{}^2 \ dx},
      \end{align*}
      for all $z \in H^2(\Omega)$ such that $\partial_{\bf n} z \vert_{\Gamma} = 0$ and $z \geq \alpha$ in $\Omega$, for some $\alpha > 0$.
    \end{lemma}

    To make the aforementioned formal calculations in a rigorous manner, we must adapt \eqref{eq_formal_u_w} to the truncated models \eqref{problema_P_m_intro}. To avoid divisions by zero, we consider the change of variables $z_m(t,x) = \sqrt{v_m(t,x) + \alpha}$, for some $\alpha > 0$ to be chosen latter. With this change of variables, we write the truncated problem \eqref{problema_P_m_intro} as the equivalent problem

    \begin{equation} \label{problema_P_u_m_z_m}
      \begin{array}{rl}
        \partial_t u_m - \Delta u_m & = - \nabla \cdot (T^m(u_m) \nabla (z_m)^2), \\
        \partial_t z_m - \Delta z_m - \dfrac{\norm{\nabla z_m}{}^2}{z_m} & 
        = - \dfrac12 T^m(u_m)^s z_m + \dfrac{\alpha}2 \dfrac{T^m(u_m)^s}{ z_m} ,\\
        \partial_{\bf n} u_m \vert_{\Gamma} & =  \partial_{\bf n} z_m \vert_{\Gamma} = 0 ,\\
        u_m(0) & = u^0_m, \quad z_m(0) = \sqrt{v^0_m + \alpha}.
      \end{array}
    \end{equation}

    Analogously to the formal calculations, now, to obtain an energy inequality, we must test the $z_m$-equations of \eqref{problema_P_u_m_z_m} by $- \Delta z_m$. Then Lemma \ref{lemma_termo_fonte_final} is applied to prove the following result.
    \begin{lemma}
      The solution $(u_m,z_m)$ of \eqref{problema_P_u_m_z_m}, satisfies the inequality
      \begin{equation*}
        \begin{array}{l}
          \dfrac{1}{2} \dfrac{d}{dt} \norma{\nabla z_m}{L^2(\Omega)}^2 + C_1 \Big ( \D{\int_{\Omega}}{\norm{D^2 z_m}{}^2 \ dx} + \D{\int_{\Omega}}{\frac{\norm{\nabla z_m}{}^4}{z_m^2} \ dx} \Big ) \\[12pt]
          + \dfrac{1}{2} \D{\int_{\Omega}{T^m(u_m)^s \norm{\nabla z_m}{}^2 \ dx}} \leq - \dfrac{s}{4} \D{\int_{\Omega}}{T^m(u_m)^{s-1} \nabla (z_m^2) \cdot \nabla T^m(u_m) \ dx} \\[12pt]
          + \dfrac{s}{2} \sqrt{\alpha} \D{\int_{\Omega}}{T^m(u_m)^{s-1} \norm{\nabla z_m}{} \norm{\nabla T^m(u_m)}{} \ dx} + C_2 \int_{\Omega}{\norm{\nabla z_m}{}^2 \ dx}.
        \end{array}
      \end{equation*}
      \label{lemma_tratamento_equacao_z}
    \end{lemma}
    Next, we must deal with the $u_m$-equation of \eqref{problema_P_u_m_z_m}. In \cite{ViannaGuillen2023uniform}, the cases $s = 1$, $s \in (1,2)$ and $s \geq 2$ are treated separately. This happens, in part, because the test functions involved have different properties concerning the singularity at zero. In fact, if we take the function $g'(u)$, given by \eqref{funcao_teste_formal}, we observe that: if $s = 1$ then $g'(u)$ and $g''(u)$ have a singularity at $u = 0$; if $s \in (1,2)$, then only $g''(u)$ is singular at $u = 0$; and if $s \geq 2$, then neither $g'(u)$ nor $g''(u)$ are singular.

    To give an idea of the procedure, we mention the case $s = 1$. In this case, to avoid divisions by zero, we test $u_m$-equation of \eqref{problema_P_u_m_z_m} by $ln(T^m(u_m) + 1)$. This gives us
    \begin{align*}
      & \frac{d}{dt} \int_{\Omega}{g_m(u_m) \ dx} + \int_{\Omega}{\frac{(T^m)'(u_m)}{T^m(u_m) + 1} \norm{\nabla u_m}{}^2 \ dx} \\
      & = \int_{\Omega}{\frac{T^m(u_m)}{T^m(u_m) + 1} \nabla (z_m^2) \cdot \nabla T^m(u_m) \ dx},
    \end{align*}
    where $g_m(r)$ is a primitive of $ln(T^m(r) + 1)$. After some manipulations, we obtain
    \begin{equation}
      \begin{array}{l}
        \dfrac{d}{dt} \D \int_{\Omega}{g_m(u_m) \ dx} + C \int_{\Omega}{\norm{\nabla [T^m(u_m) + 1]^{1/2}}{}^2 \ dx} \\[6pt]
        \leq \D \int_{\Omega}{\nabla (z_m^2) \cdot \nabla T^m(u_m) \ dx} + \tilde{C} \norma{\nabla z_m}{L^2(\Omega)}^2.
      \end{array}
      \label{estimativa_u_m_1_s=1}
    \end{equation}
    Summing \eqref{estimativa_u_m_1_s=1} to $1/4$ times the inequality of Lemma \ref{lemma_tratamento_equacao_z}, for $s = 1$, the terms related to $\int_{\Omega}{\nabla (z_m^2) \cdot \nabla T^m(u_m) \ dx}$ cancel each other. Estimating the terms which do not cancel, the following energy inequality is proved.
    \begin{lemma}[\bf Energy inequality for $\boldsymbol{s = 1}$]
        The solution $(u_m, z_m)$ of the problem \eqref{problema_P_u_m_z_m} satisfies, for sufficiently small $\alpha > 0$,
        \begin{equation}
          \begin{array}{c}
            \dfrac{d}{dt} \left[ \dfrac{1}{4} \D{\int_{\Omega}}{g_m(u_m) \ dx} + \dfrac{1}{2} \D{\int_{\Omega}}{\norm{\nabla z_m}{}^2 \ dx} \right] \\[12pt]
            + C \D{\int_{\Omega}}{\norm{\nabla [T^m(u_m) + 1]^{1/2}}{}^2 \ dx} + \frac{1}{4} \D{\int_{\Omega}}{T^m(u_m) \norm{\nabla z_m}{}^2 \ dx} \\[12pt]
            + C_1 \Big ( \D{\int_{\Omega}}{\norm{D^2 z_m}{}^2 \ dx} + \D{\int_{\Omega}}{\frac{\norm{\nabla z_m}{}^4}{z_m^2} \ dx} \Big ) \leq C \int_{\Omega}{\norm{\nabla z_m}{}^2 \ dx}.
          \end{array}
          \label{estimativa_u_v_m_1_s=1}
        \end{equation}
        \label{lemma_estimativa_u_v_m_1_s=1}
      \end{lemma}
      From \eqref{estimativa_u_v_m_1_s=1} and the corresponding energy inequalities for the cases $s \in (1,2)$ and $s \geq 2$, $m$-independent estimates for $T^m(u_m)$ and $\nabla z_m$ are obtained in \cite{ViannaGuillen2023uniform}, which also us gives estimates for $\nabla v_m$. With these new estimates, we go back to the $u_m$-equation of \eqref{problema_P_m_intro} and get $m$-independent estimates for $u_m$.

      Using all the aforementioned $m$-independent estimates, it is possible to pass to the limit in \eqref{problema_P_m_intro}, proving Theorem \ref{theorem_3D_existencia}. To prove Theorem \ref{theorem_2D_existencia_unicidade_estimativas}, we obtain stronger $m$-independent estimates for $\{ u_m \}_m$ and $\{ v_m \}_m$ which are available only in $2$D domains, obtaining stronger regularity for the solution of \eqref{problema_P}. Finally, with this stronger regularity, uniqueness of solution is proved.


\section{Discrete Schemes}
\label{sec:TimeDiscrete}

  The numerical approximation of chemotaxis models is a relevant and growing research topic (see \cite{marrocco2003numerical,saito2007conservative,epshteyn2009fully,saito2011error,ibrahim2014efficacy,bessemoulin2014finite,zhang2016characteristic,chertock2018high,guillen2019unconditionally,guillen2020study2,guillen2022comparison,guillen2020theoretical,guillen2021chemorepulsion,gutierrez2021analysis,badia2022bound,acosta2023unconditionally} and the references cited therein). Nevertheless, when we turn to problem \eqref{problema_P}, we still find a relatively small amount of studies on its numerical approximation. To the best of our knowledge, we can cite \cite{duarte2021numerical,guillen2023finite} about the numerical approximation of \eqref{problema_P}, both just for the case $s = 1$.
  
  In \cite{duarte2021numerical} a chemotaxis-Navier-Stokes system is approached via Finite Elements (FE). In fact, by assuming the existence of a sufficiently regular solution, if the initial data of the scheme are small perturbations of the initial data of this regular solution, optimal error estimates are deduced. Continuing with this analysis of approximating regular solutions, in \cite{feng2021error} a decoupled linear positivity preserving FEM scheme for a chemotaxis-Stokes problem is proposed. In \cite{li2023error}, the authors claim to have enhanced \cite{duarte2021numerical} in the sense that they avoid introducing an auxiliary variable in the design of their scheme. In \cite{beltran2023chemotaxis} a more general chemotaxis-Navier-Stokes system with Lotka-Volterra competitive term is studied. The authors establish an additional regularity hypothesis that allows them to prove that a weak solution is actually a strong solution and propose a fully discrete scheme to approximate this more regular solution.
    
  This idea of approximating exact solutions that are more regular than the weak solutions has also been used previously in \cite{epshteyn2009fully,zhang2016characteristic}, for example, where the authors assume boundedness of the exact solution of a Keller-Segel chemotaxis model. The drawback of this kind of results is that the existence of such a regular solution is not clear in general, especially when we consider polyhedral domains, which are broadly used in numerical simulations. Besides, the design of a discrete scheme converging towards a weak solution of \eqref{problema_P} is not addressed by these works.
    
  In \cite{guillen2023finite}, motivated by the treatment given to the chemorepulsion model with linear production in \cite{guillen2019unconditionally}, some FE schemes are designed to approximate \eqref{problema_P}, with $s = 1$. The authors focus on FE schemes satisfying properties such as conservation of cells, discrete energy law and approximate positivity rather than convergence. There is evidence that the preservation of such properties, especially the positivity, could possibly enhance the performance of the numerical schemes, avoiding spurious oscillations, even if the numerical scheme is convergent, which is the case in \cite{guillen2022comparison}. Numerical simulations are carried out to compare the performance of the different schemes. 
  In particular, a scheme satisfying a discrete energy law that, in 1D domains, lead to decreasing energy is presented in \cite{guillen2023finite}. The convergence of the presented schemes in $2$D or $3$D domains, however, is not clear.

  One of the main difficulties of addressing issues concerning the convergence of numerical schemes towards weak solutions of \eqref{problema_P} is probably the lack of energy (\emph{a priori}) estimates over the solutions of the schemes. Even if we consider only time discretizations of \eqref{problema_P}, the task of designing a convergent scheme is not straightforward. This could be attributed to the complex technique needed in order to obtain energy estimates through the cancellation of the chemoattraction and consumption effects.


  To the best of our knowledge, excepting the case of $1$D domains \cite{guillen2023finite}, the only time discrete scheme for \eqref{problema_P} possessing an energy inequality from which one can obtain convergence to weak solutions is proposed in \cite{guillen2023convergence}. This convergence is valid in $3$D domains and is based on energy estimates. Moreover, the proposed scheme preserves the properties of positivity and conservation of the population of cells. 
  We finish this section by showing some results of \cite{guillen2023convergence} in more detail to highlight the procedures and difficulties involved in the design of convergent discrete schemes approximating the chemotaxis-consumption models discussed here.

  The first issue that we point out is that the design of the time discrete scheme is based on the analysis that was carried out in \cite{ViannaGuillen2023uniform}. According to Section \ref{sec:Analytical}, in \cite{ViannaGuillen2023uniform}, it was convenient to rewrite \eqref{problema_P} in terms of the variable $z = \sqrt{v+\alpha^2}$, because the test functions involved in obtaining a discrete energy law become simpler. To adapt this procedure to a time discrete scheme, it was found to be more appropriate to propose a scheme using the variables $(u,z)$, instead of $(u,v)$.
    
  In great part, it is due to the following fact. To discretize in time, we will divide the interval $[0,\infty)$ in subintervals denoted by $I_n = (t_{n-1}, t_n)$, with $t_0 = 0$ and $t_n = t_{n-1} + k$, where $k > 0$ is the time step. If $\{ z^n \}_n$ is a sequence of functions defined in $\Omega$, then we use the following notation for the discrete time derivative
  \begin{equation}  
    \delta_t z^n = \frac{z^n - z^{n-1}}{k}, \quad \forall n \geq 1 .
  \end{equation} 
  Then, the following lemma is usually applied to deal with the discrete time derivative term.
  \begin{lemma}[\bf \cite{eyre1998unconditionally}] \label{lema_delta_t}
    Let $z^n, z^{n-1} \in L^\infty(\Omega)$ and let $F: \mathbb{R} \rightarrow \mathbb{R}$ be a $C^2$ function. Then
    \begin{equation*}
      \delta_t z^n \, F'(z^n) = \delta_t F(z^n) + \frac{1}{2k} F''(c^n(x))(z^n(x) - z^{n-1}(x))^2,
    \end{equation*}  
    where $c^n(x)$ is an intermediate value between $z^n(x)$ and $z^{n-1}(x)$.
  \end{lemma}
  Hence, if a time discrete scheme in the variables $(u,v)$ is proposed and we want to adapt the procedures used in \cite{ViannaGuillen2023uniform} to obtain estimates for the discrete solutions, accounting for what has been shown in Section \ref{sec:Analytical}, formally, we would have to start by multiplying the $v^n$-equation of the time discrete scheme by $F'(v^n) = 1/\sqrt{v^n}$ (assuming $v^j(x) > 0$, for all $x \in \Omega$ and for all $j$). By Lemma~\ref{lema_delta_t}, we would find the following:
  \begin{equation*}
    \delta_t v^n \, \frac{1}{\sqrt{v^n(x)}} = 2\, \delta_t (\sqrt{v^n(x)})
     - \frac{1}{4k} \frac{1}{\sqrt{c^n(x)}}(v^n(x) - v^{n-1}(x))^2.
  \end{equation*}
  This means that we would have to treat the negative term
  \begin{equation*}
    - \frac{1}{4k} \frac{1}{\sqrt{c^n(x)}}(v^n(x) - v^{n-1}(x))^2.
  \end{equation*}
  Since it is not clear how this term could be estimated, we circumvent this difficulty by proposing the time discrete scheme using directly the variable $(u,z)$.   
  Therefore, the following reformulation of \eqref{problema_P} is considered in \cite{guillen2023convergence},
  {\small \begin{equation}  
    \left\{\begin{array}{l}
      \partial_t u - \Delta u  = -  \nabla \cdot (u \nabla (z)^2), \quad
      \partial_t z - \dfrac{\norm{\nabla z}{}^2}{z} - \Delta z  = - \dfrac{1}{2} u^s 
      \left(z - \dfrac{\alpha^2}{z} \right) , \\
      \partial_\eta u \Big \vert_{\Gamma}  =  \partial_\eta z \Big \vert_{\Gamma} = 0, \quad
      u(0)  = u^0, \quad z(0) = \sqrt{v^0 + \alpha^2},
    \end{array}\right.
    \label{problema_P_u_z_alpha} 
  \end{equation}} 
where $\alpha > 0$ is a fixed real number to be chosen later. Since it is proved in \cite{ViannaGuillen2023uniform} that the $v$-equation of \eqref{problema_P} is satisfied pointwisely, with $v \in L^2(0,\infty;H^2(\Omega))$, one can check by straightforward calculations that \eqref{problema_P_u_z_alpha} is equivalent to \eqref{problema_P} if we use the change of variables $z = \sqrt{v + \alpha^2}$. We summarize this statement in the following lemma for further use.
    
  \begin{lemma} \label{lema_equivalencia_problemas}
    Problems \eqref{problema_P} and \eqref{problema_P_u_z_alpha} are equivalent. More precisely, $(u,z)$ is a weak solution of \eqref{problema_P_u_z_alpha} if, and only if, $(u,v)$ is a weak solution of \eqref{problema_P}, with $v = z^2 - \alpha^2$.
  \end{lemma}
    
  We will also use the following upper truncation of $u$, for each fixed $m > 0$, 
  \begin{equation*} 
    T^m(u) = \left \{
    \begin{array}{cl}
      u, & \mbox{ if } u \leq m, \\
      m, & \mbox{ if } u \geq m.
    \end{array}
    \right .
  \end{equation*} 
    
  Then, we propose the following time discrete scheme:

\
    
  \noindent  {\bf Initialization: } Consider the initial conditions $u^0 \in L^2(\Omega)$, $z^0 = \sqrt{v^0 + \alpha^2} \in L^{\infty}(\Omega)$ with $v^0 \in L^{\infty}(\Omega)$.
  
  \noindent  {\bf Step $n$: }    
  Given non-negative functions $u^{n-1} \in L^2(\Omega)$, $z^{n-1} \in L^{\infty}(\Omega)$ and $v^{n-1} \in L^{\infty}(\Omega)$, 
  \begin{enumerate}
    \item Find  $(u^n,z^n) \in H^2(\Omega)^2$, satisfying the bounds 
    $$u^n(x) \geq 0\quad \hbox{and}\quad \norma{z^{n-1}}{L^{\infty}(\Omega)} \geq z^n(x) \geq \alpha \quad \hbox{a.e. $x \in \Omega$,}
    $$
    and the boundary-value problem 
    \begin{equation}  \label{NLTD}
      \left\{
      \begin{array}{l}
        \delta_t u^n - \Delta u^n =  \nabla \cdot \Big (  T^m(u^n) \nabla (z^n)^2 \Big ), \\
        \delta_t z^n - \dfrac{\norm{\nabla z^n}{}^2}{z^n} - \Delta z^n = - \dfrac{1}{2} T^m(u^n)^s \left(z^n - \dfrac{\alpha^2}{z^n} \right), \\
        \partial_{\eta} u^n \Big \vert_{\partial \Omega} = \partial_{\eta} z^n \Big \vert_{\partial \Omega} = 0.
      \end{array}
      \right.
    \end{equation} 
    \item Two variants for the approximation of $v$ are possible (equally denoted), either depending on $z^n$ or $u^n$:
    \begin{itemize}
      \item  Find $v^n = v^n(z^n) \in H^2(\Omega)$ as
      \begin{equation}  \label{v-z}
        v^n = (z^n)^2 - \alpha^2 .
      \end{equation} 
      \item Find $v^n = v^n(u^n) \in H^2(\Omega)$ as the unique solution of the linear problem
      \begin{equation}  \label{v-u}
        \delta_t v^n - \Delta v^n + T^m(u^n)^s v^n = 0,
        \quad 
        \partial_{\eta} v^n \Big \vert_{\partial \Omega} =0.
      \end{equation} 
    \end{itemize}
  \end{enumerate}

\

  \begin{remark} \label{remark_dependencia_esquema_em_m}
    Note that $u^n$, $z^n$ and $v^n$ depend on the truncation parameter $m$. For simplicity, from now on, we consider that $m \in \mathbb{N}$.
    \hfill $\square$
  \end{remark}

\

  Once defined the time discrete scheme, the first step to its analysis is the existence of discrete solutions.
  \begin{theorem}{\bf (\cite{guillen2023convergence} Existence of solution of \eqref{NLTD})}
    Suppose $(u^{n-1},z^{n-1}) \in L^2(\Omega) \times L^{\infty}(\Omega)$ with $u^{n-1}(x) \geq 0$ and $z^{n-1}(x) \geq \alpha$ $a.e.$ $x \in \Omega$. Then there is a solution $(u^n,z^n)$ of \eqref{NLTD} which satisfies $u^n(x) \geq 0$ and $\norma{z^{n-1}}{L^{\infty}(\Omega)} \geq z^n(x) \geq \alpha$ $a.e.$ $x \in \Omega$.
    \label{teo_existencia_NLTD}
  \end{theorem}

  Using the functions $u^n$, $z^n$ and $v^n$ introduced above in step $n$, we define the piecewise function $u_m^k$ and the locally linear and globally continuous function $\tilde{u}_m^k$ by
  \begin{equation}  \label{funcao_u^kr_u^k}
    \begin{array}{c}
      u_m^k(t,x) = u^n(x) \mbox{ and } \\[6pt]
      \tilde{u}_m^k(t,x) = u^n(x) + \dfrac{(t - t_n)}{k} \big ( u^n(x) - u^{n-1}(x) \big ), \mbox{ if } t \in [t_{n-1},t_n).
    \end{array}
  \end{equation} 
  Analogously, we define the functions $z_m^k$, $\tilde{z}_m^k$, $v_m^k$ and $\tilde{v}_m^k$.
  
  Now, we are in position to present the main result that is proved along \cite{guillen2023convergence}.
  \begin{theorem} \label{teo_principal}
    For each $n \in \mathbb{N}$, there exists at least one solution $(u^n,z^n)$  of \eqref{NLTD}, that jointly to $v^n$ defined by \eqref{v-z} or \eqref{v-u} leads us to $(u^n,v^n)$ satisfying
    $$u^n(x) \geq 0, \qquad \norma{v^0}{L^{\infty}(\Omega)} \geq v^n(x) \geq 0 \quad a.e. \ x \in \Omega.$$
    Moreover, up to a subsequence, $(u_m^k,v_m^k)$ converges towards a weak solution $(u,v)$ of \eqref{problema_P} as $(m,k) \to (\infty,0)$.
  \end{theorem}

\

  \begin{remark}
    We have some points to remark:
    
    \begin{enumerate}
      \item The number $\alpha$ is a sufficiently small positive real number that is chosen in Lemmas \ref{lemma_estimativa_u_v_m_1_s_intermediario_time_discrete} and \ref{lemma_estimativa_u_v_m_1_s_geq_2_time_discrete} independently of $m$ and $k$.
      \item The convergence result given in Theorem \ref{teo_principal} as $(m,k) \to (\infty,0)$ is unconditional, that is, there is not any constraint over $m$ and $k$ as long as $m \to \infty$ and $k \to 0$.
      \item In particular,  we also prove the  result on existence of weak solutions to \eqref{problema_P} in $3$D domains given in \cite[Theorem 1]{ViannaGuillen2023uniform} but, this time, as a consequence of the convergence of the time discrete scheme introduced in this paper.
      \item In $2$D domains, there exists a unique strong solution of \eqref{problema_P} (see \cite{ViannaGuillen2023uniform}). The proof is achieved through the obtaining of stronger $m$-independent estimates for the solution of an adequate truncated problem. Unfortunately, it is not clear how we could adapt these strong estimates for the time discrete scheme. Consequently, in $2$D, the convergence of the whole sequence of solutions of the time discrete scheme \eqref{NLTD} towards the unique strong solution of \eqref{problema_P} as $(m,k) \to (0,\infty)$ remains as an open problem.
    \end{enumerate}
    \hfill $\square$
  \end{remark}


  We highlight the main ideas of the proof of Theorem \ref{teo_principal}. We first use the global in time functions introduced in \eqref{funcao_u^kr_u^k} to rewrite \eqref{NLTD} as the following differential system, $a.e.$ in $(t,x) \in (0,\infty) \times \Omega$,
  \begin{equation}  
    \left\{
    \begin{array}{l}
      \partial_t \tilde{u}_m^k - \Delta u_m^k  =  \nabla \cdot \Big (  T^m(u_m^k) \nabla (z_m^k)^2 \Big ), \\
      \partial_t \tilde{z}_m^k - \dfrac{\norm{\nabla z_m^k}{}^2}{z_m^k} - \Delta z_m^k  = - \dfrac{1}{2} T^m(u_m^k)^s \left(z_m^k - \dfrac{\alpha^2}{z_m^k}\right).
    \end{array}
    \right.
    \label{NLTD_equiv}
  \end{equation}
  Then, we must pass to the limit in \eqref{NLTD_equiv} as $(m,k) \to (\infty,0)$. To do so, we need adequate $(m,k)$-independent bounds for the sequences $\{ u_m^k \}$, $\{ z_m^k \}$, $\{ \tilde{u}_m^k \}$ and $\{ \tilde{z}_m^k \}$, which will be obtained from estimates involving $(u^n,z^n)$, for each $n \in \mathbb{N}$.
  
  Analogously to Section \ref{sec:Analytical}, we have some direct $(m,k)$-independent estimates for the approximate solutions which are summarized in the following result.
  \begin{lemma}{\bf ($\boldsymbol{(m,k,n)}$-uniform estimates)} \label{lemma_limitacao_grad_z_uniforme_wrt_k_m_T}
    Let $(u^n,z^n)$ be a solution of \eqref{NLTD}. Then one has
    \begin{enumerate}
      \item $\D{\int_{\Omega}} u^n \ dx = \D{\int_{\Omega}} u^0  dx$, \ for all $n \in \mathbb{N}$;
      \item $\norma{z^n}{L^2(\Omega)}^2 + \D{\sum_{j=1}^n} \norma{z^j - z^{j-1}}{L^2(\Omega)}^2 \leq \norma{z^0}{L^2(\Omega)}^2$, \ for all $n \in \mathbb{N}$;
      \item $k \D{\sum_{j=1}^n}{\norma{\nabla z^j}{L^2(\Omega)}^2} \leq \frac{1}{4 \alpha^2} \norma{v^0 + \alpha^2}{L^2(\Omega)}^2$, \ for all $n \in \mathbb{N}$.
    \end{enumerate}
  \end{lemma}
  Using the above estimates for $(u^n,v^n)$ we conclude the following $(m,k)$-independent bounds for $u_m^k$, $z_m^k$, $\tilde{u}_m^k$ and $\tilde{z}_m^k$:
  \begin{equation*}
    u_m^k, \tilde{u}_m^k \mbox{ are bounded in } L^{\infty}(0,\infty;L^1(\Omega)),
  \end{equation*}
  \begin{equation*}
    z_m^k, \tilde{z}_m^k \mbox{ are bounded in } L^{\infty}(0,\infty;L^2(\Omega)),
  \end{equation*}
  \begin{equation*}
    \nabla z_m^k, \nabla \tilde{z}_m^k \mbox{ are bounded in } L^2(0,\infty;L^2(\Omega)).
  \end{equation*}

  These bounds are not sufficient to pass to the limit as $(m,k) \to (\infty,0)$ and stronger estimates involving $(u^n,z^n)$ are necessary. The procedures applied to obtain such estimates are  based on the procedures applied to the continuous problem shown in Section \ref{sec:Analytical}.

  First, a time discrete version of Lemma \ref{lemma_tratamento_equacao_z} is proved in \cite{guillen2023convergence}. In fact, by testing the $z^n$-equation of \eqref{NLTD} by $- \Delta z^n$, applying Lemma \ref{lemma_termo_fonte_final} and performing a series of appropriate manipulations, one can prove the following lemma.
  \begin{lemma} \label{lemma_tratamento_equacao_z_n}
    Any solution $(u^n,z^n)$ of \eqref{NLTD}, satisfies the inequality
    \begin{equation}  
      \begin{array}{l}
        \dfrac{1}{2} \delta_t \norma{\nabla z^n}{L^2(\Omega)}^2 + \dfrac{1}{2 k} \norma{\nabla z^n - \nabla z^{n-1}}{L^2(\Omega)}^2 \\
        + C_1 \Big ( \D{\int_{\Omega}}{\norm{D^2 z^n}{}^2 \ dx} + \D{\int_{\Omega}}{\frac{\norm{\nabla z^n}{}^4}{(z^n)^2} \ dx} \Big ) + \dfrac{1}{2} \D{\int_{\Omega}{T^m(u^n)^s \norm{\nabla z^n}{}^2 \ dx}} \\
        \leq \dfrac{s}{4} \D{\int_{\Omega}}{T^m(u^n)^{s-1} \nabla (z^n)^2 \cdot \nabla T^m(u^n) \ dx} \\
        + \dfrac{s \alpha}{2} \D{\int_{\Omega}}{T^m(u^n)^{s-1} \norm{\nabla z^n}{} \norm{\nabla T^m(u^n)}{} \ dx} + C_2 \D{\int_{\Omega}}{\norm{\nabla z^n}{}^2 \ dx}.
      \end{array}
      \label{eq_tratamento_equacao_z}
    \end{equation}
  \end{lemma}

  Next, a local energy inequality for $(u^n, z^n)$ based on the cancellation between the chemotaxis and consumption effects is obtained. Here, the analysis must be divided into the cases $s \in [1,2)$ and $s \geq 2$. We consider the function $g_m$ defined by $g_m(r) = \D{\int_0^r}{g'_m(\theta) \ d\theta}$, where
  \begin{equation*}
    g_m'(r) = \left \{ 
    \begin{array}{rl}
      ln(T^m(r)), & \mbox{if } s = 1, \mbox{ for } r > 0, \\
        \dfrac{T^m(r)^{s-1}}{(s-1)}, & \mbox{if } s > 1. 
    \end{array}
    \right.
  \end{equation*} 

  By testing the $u^n$-equation of \eqref{NLTD} by $g_m'(u^n)$ and combining the resulting expression with \eqref{eq_tratamento_equacao_z}, one obtains the following energy inequalities.
  
  \begin{lemma}[\bf Energy inequality for $\boldsymbol{s \in [1,2)}$] \label{lemma_estimativa_u_v_m_1_s_intermediario_time_discrete}
    Any solution $(u^n, z^n)$ of the problem \eqref{NLTD} satisfies, for sufficiently small $\alpha^2 > 0$,
    \begin{equation}  
      \begin{array}{l}
        \delta_t \Big [ \dfrac{s}{4} \D{\int_{\Omega}}{g_m(u^n) \ dx} + \frac{1}{2} \norma{\nabla z^n}{L^2(\Omega)}^2 \Big ] \\[6pt]
        + \dfrac{1}{2 k} \norma{\nabla z^n - \nabla z^{n-1}}{L^2(\Omega)}^2 + \dfrac{1}{4} \D{\int_{\Omega}}{T^m(u^n)^s \norm{\nabla z^n}{}^2 \ dx} \\[6pt]
        + C_1 \Big ( \D{\int_{\Omega}}{\norm{D^2 z^n}{}^2 \ dx} + \D{\int_{\Omega}}{\frac{\norm{\nabla z^n}{}^4}{(z^n)^2} \ dx} \Big ) \leq C \norma{\nabla z^n}{L^2(\Omega)}^2.
      \end{array}
      \label{estimativa_u_v_m_1_s_intermediario_time_discrete}
    \end{equation} 
  \end{lemma}

  \begin{lemma}[\bf Energy inequality for $\boldsymbol{s \geq 2}$] \label{lemma_estimativa_u_v_m_1_s_geq_2_time_discrete}
    Any solution $(u^n, z^n)$ of the problem \eqref{NLTD} satisfies
    \begin{equation} \vspace{-1pt}  \label{estimativa_u_v_m_1_s_geq_2_time_discrete}
      \begin{array}{l}
        \delta_t \Big [ \dfrac{s}{4} \D{\int_{\Omega}} g_m(u^n) \ dx  + \frac{1}{2} \norma{\nabla z^n}{L^2(\Omega)}^2 \Big ] + \dfrac{1}{2 k} \norma{\nabla z^n - \nabla z^{n-1}}{L^2(\Omega)}^2 \\[6pt] + \D{\int_{\Omega}}{\norm{\nabla [T^m(u^n)]^{s/2}}{}^2  dx} + \dfrac{1}{4} \D{\int_{\Omega}}{T^m(u^n)^s \norm{\nabla z^n}{}^2 \ dx} \\[6pt]
        + C_1 \Big ( \D{\int_{\Omega}}{\norm{D^2 z^n}{}^2 \ dx} + \D{\int_{\Omega}}{\frac{\norm{\nabla z^n}{}^4}{(z^n)^2} \ dx} \Big ) \leq C \norma{\nabla z^n}{L^2(\Omega)}^2.
      \end{array}
    \end{equation} 
  \end{lemma}

  Summing the energy inequalities  \eqref{estimativa_u_v_m_1_s_intermediario_time_discrete} and \eqref{estimativa_u_v_m_1_s_geq_2_time_discrete} from $1$ to $n$, for all $n \in \mathbb{N}$, 
  and using that, by Lemma \ref{lemma_limitacao_grad_z_uniforme_wrt_k_m_T},
  \begin{equation*}
    k \sum_{n=1}^{\infty}{\norma{\nabla z^n}{L^2(\Omega)}^2}
    \quad 
    \hbox{is bounded,}
  \end{equation*}
   one obtains $(m,k)$-independent bounds for the sequences $\{ T^m(u_m^k) \}$, $\{ \nabla z_m^k \}$, $\{ T^m(\tilde{u}_m^k) \}$ and $\{ \nabla \tilde{z}_m^k \}$, which will also be uniform in time.
  To obtain $(m,k)$-independent bounds for $\{ u_m^k \}$ and $\{ \tilde{u}_m^k \}$, one has to go back to the $u^n$-equation of \eqref{NLTD} and test it by $u^n$, if $s \geq 2$, or $g'(u^n + 1)$, if $s \in [1,2)$, where
  \begin{equation*}
    g'(r) =  \left \{
    \begin{array}{cc}
      ln(r) & \mbox{ if } s = 1, \ \forall r > 0, \\
      r^{s-1}/(s-1) & \mbox{ if } s \in (1,2).
    \end{array}
    \right .
  \end{equation*}

  Analogously to Section \ref{sec:Analytical}, the aforementioned $(m,k)$-independent bounds are used in order to pass to the limit in \eqref{NLTD_equiv}. However, now there is an important difference. We are capable of obtaining strong convergence for the sequences of functions $\tilde{u}^k_m$ and $\tilde{z}^k_m$, whose time derivatives can be bounded independently of $(m,k)$ by going back to \eqref{NLTD_equiv}, but we need strong convergence of the functions $u^k_m$ and $z^k_m$, which appear in the nonlinear terms. A usual way to deal with this issue is to prove that $u^k_m - \tilde{u}^k_m$ and $z^k_m - \tilde{z}^k_m$ converge to zero in adequate norms as $(m,k) \to (\infty,0)$. 
  
  Summing \eqref{estimativa_u_v_m_1_s_intermediario_time_discrete} from $1$ to $\infty$, we conclude, in particular, that
  \begin{equation}  \label{limitacao_diferenca_z_H_1}
    \D{\sum_{j=1}^{\infty}}{\norma{z^j - z^{j-1}}{H^1(\Omega)}^2} \leq C.
  \end{equation} 
  Then, using \eqref{limitacao_diferenca_z_H_1}, there is a positive constant $C$, independent of $m$ and $k$, such that
  \begin{equation}  \label{convergencia_diferenca_z_m_em_H_1}
    \norma{z_m^k - \tilde{z}_m^k}{L^2(0,\infty;H^1(\Omega))}^2 \leq C\, k,
  \end{equation}
  that is, $\norma{z_m^k - \tilde{z}_m^k}{L^2(0,\infty;H^1(\Omega))} \to 0$ as $(m,k) \to (\infty,0)$. Accounting for $u^k_m - \tilde{u}^k_m$, in case $s \geq 2$, the result is also straightforward. As it is shown in \cite{guillen2023convergence}, from the estimate
  \begin{equation*} 
    \D{\sum_{j=1}^{\infty}}{\norma{u^j - u^{j-1}}{L^2(\Omega)}^2} \leq C,
  \end{equation*}  
  one has that
  \begin{equation} \label{convergencia_diferenca_u_m_em_L_2} 
    \norma{u_m^k - \tilde{u}_m^k}{L^2(0,\infty;L^2(\Omega))}^2 \leq C k.
  \end{equation}

  But, for the case $s \in [1,2)$, the estimate we have to work with is
  \begin{equation*} 
    \sum_{j=1}^{\infty} \int_{\Omega} \frac{1}{(c^j(x) + 1)^{2 - s}} (u^j(x) - u^{j-1}(x))^2 dx \leq C,
  \end{equation*}
  where $c^j(x)$ is a point between $u^{j-1}(x)$ and $u^j(x)$. The weight $(c^j(x) + 1)^{s-2}$ can degenerate and then a conclusion such as \eqref{convergencia_diferenca_z_m_em_H_1} and \eqref{convergencia_diferenca_u_m_em_L_2} is not clear. In this case the analysis is more delicate and the result that is proved in \cite{guillen2023convergence} is the following.
  \begin{lemma}{\bf (\cite{guillen2023convergence})} \label{lema_f''}
    There is a positive constant $C$, independent of $m$ and $k$, such that
    \begin{equation*} 
      \norma{u_m^k - \tilde{u}_m^k}{L^2(0,\infty;L^s(\Omega))}^2 \leq C\, k.
    \end{equation*} 
  \end{lemma}
  Accounting for all that has been exposed so far, in \cite{guillen2023convergence} it is proved, for any $s \geq 1$ fixed, that, up to a subsequence, $(u_m^k, z_m^k)$ converges to a weak solution $(u,z)$ of \eqref{problema_P_u_z_alpha} as $(m,k) \to (\infty, 0)$.
  
  Finally, to conclude the proof of Theorem \ref{teo_principal}, it is shown in \cite{guillen2023convergence} that $(u_m^k, v_m^k)$ converges to a weak solution $(u,v)$ of \eqref{problema_P} as $(m,k) \to (\infty, 0)$, where $v_m^k$ is given either by  \eqref{v-z} or by \eqref{v-u}.


\section{Optimal Control} \label{sec:Control}

  The works cited so far address the analysis of the chemotaxis-consumption model \eqref{problema_P} and help us understand how the system evolves over time from the given initial data, without  external interference. However, when it comes to PDE models describing physical phenomena, as important as the analysis of the system itself, are the studies of control problems related to them.
  
  Concerning chemotaxis models, particularly relevant are the optimal control problems. Due to the low regularity of controlled problems and the lack of uniqueness in $3$D of the control-to-state operator, some of the existing works are concentrated on optimal control problems in $2$D domains, where one  usually has existence and uniqueness of a strong solution to the controlled model, which allows to show the existence of global optimal solution and to derive an optimality system, establishing existence and regularity of Lagrange multipliers associated to any local optimum.
  
  For more details on this kind of works, we refer the interested reader to the works on control problems in $2$D domains related to: a Keller-Segel model \cite{ryu2001optimal}; a chemorepulsion-production model \cite{guillen2020optimal,guillen2020bi}; a Keller-Segel logistic model \cite{silva2022bilinear}; a chemotaxis model with indirect consumption \cite{yuan2022optimal}; and a chemotaxis-haptotaxis model \cite{tang2022optimal}.
  
  When we turn to optimal control problems related to chemotaxis models in $3$D domains, this analysis becomes more complex. In great part, this is because, in $3$D domains, we have existence of weak solutions, however, in many cases, there is not any result on the existence and uniqueness of global in time strong solutions. To overcome this difficulty some authors introduce a regularity criterion, which is a mild additional regularity hypothesis on a weak solution, sufficient to conclude that this weak solution is actually the unique strong solution.
  
  For a motivated introduction of this technique we refer the reader to \cite{casas1998optimal}, where the author studied an optimal control problem related to the Navier-Stokes equations in $3$D domains. For chemotaxis related works in this setting, in which a regularity criterion is established to study the optimal control problem in $3$D domains, we refer the reader to: \cite{guillen2020regularity}, for a chemorepulsion-production model; \cite{lopez2021optimal}, for a chemotaxis-Navier-Stokes model; and \cite{guillen2023optimal}, for the chemotaxis-consumption model \eqref{problema_P_controlado}.

  The drawback of using a regularity criterion is that it is not clear if the admissible set is nonempty, except for some specific cases. In \cite{guillen2020regularity}, the authors show that if $\Omega_c = \Omega$, that is, if the control acts in the whole domain, and the initial chemical density $v^0$ is strictly positive and separated from zero, then the admissible set is nonempty. To do that, the idea is to define the control $f$ \emph{a posteriori}, depending on a regular pair $(u,v)$. In \cite{lopez2021optimal} and \cite{guillen2023optimal}, to prove that the admissible set is nonempty, it is sufficient that $\Omega_c = \Omega$ (without imposing strict positivity for $v^0$).

  Accounting for that, in \cite{guillen2022optimal}, an optimal control problem related to \eqref{problema_P_controlado} is studied in the weak solution setting, that is, without using any regularity criterion or hypothesis over the admissible set. First, the concept of weak solutions of the controlled model \eqref{problema_P_controlado} satisfying an energy inequality is considered. Next, an optimal control problem for which it is possible to prove existence of global optimal solution is considered and, to conclude, the relation between this optimal control problem and two other related ones, where the existence of optimal solution can not be proved, is discussed.

  We remark that, in this weak solution setting, it is not clear how to obtain first order optimality conditions. In fact, since we only have  weak regularity, it is not possible to prove the well-posedness of the linearized problem around a local optimal solution, which is a fundamental step to apply a generic Lagrange multiplier method as in \cite{guillen2023optimal}. Also, given a control $f$, the uniqueness of the state $(u,v)$ can not be guaranteed, hence it is not possible to define the  ``control-to-state'' mapping and follow the procedure used in \cite{lopez2021optimal} to compute the derivative of the state with respect to the control.
  
  Since we are focused on the chemotaxis-consumption model \eqref{problema_P_controlado}, in subsections \ref{subsection:control weak} and \ref{subsection:control strong} we will see in more detail the results of \cite{guillen2022optimal} and \cite{guillen2023optimal}, respectively. First, it is convenient to restate here the studied controlled model, in $Q = (0,T) \times \Omega$,
  \begin{equation} \label{problema_P_controlado_2}
    \left\{\begin{array}{l}
      \partial_t u - \Delta u  = - \nabla \cdot (u \nabla v), \quad
      \partial_t v - \Delta v  = - u^s v + f v 1_{\Omega_c}, \\ [6pt]
      \partial_{\bf n} u  \vert_{\Gamma}  =  \partial_{\bf n} v  \vert_{\Gamma} = 0, \quad
      u(0)  = u^0, \quad v(0) = v^0,
    \end{array}\right.
  \end{equation}
  and to introduce some important necessary assumptions,
  \begin{equation} \label{assumptions_set}
    \begin{array}{c}
      \Omega \subset \mathbb{R}^3 \mbox{ is a bounded domain with boundary } \Gamma \mbox{ of class } C^{2,1}, \\[6pt]
      \Omega_c \subset \Omega \mbox{ is a subdomain with boundary } \Gamma_c \mbox{ locally Lipschitz}, \\[6pt]
      s \geq 1 \mbox{ and } f \in L^q(Q_c) \mbox{ with }  q > 5/2,
    \end{array}
  \end{equation}
  as well as to clarify the concepts of weak and strong solution in the following definitions. First, let us define some functional spaces related to the concept of weak solutions. For $s \in [1,2)$,
  \begin{equation*}
    \begin{array}{rl}
      \mathcal{U} = & \hspace{-4pt} \Big \{ u \ \vert \ u \in L^{\infty}(0,T;L^s(\Omega)) \cap  L^{5s/3}(Q), \  \nabla u \in   L^{5s/(3 + s)}(Q), \\
      & \qquad \partial_t u \in L^{5s/(3 + s)}(0,T;W^{1,5s/(4s-3)}(\Omega)') \Big \},
    \end{array}
  \end{equation*}
  for $s \geq 2$,
  \begin{equation*}
    \begin{array}{rl}
      \mathcal{U} = & \Big \{ u\ \vert \ u \in L^{\infty}(0,T;L^s(\Omega)) \cap L^{5s/3}(Q),\  \nabla u \in L^{2}(Q), \\
      & \qquad  \partial_t u \in L^2(0,T;H^1(\Omega)') \Big \},
    \end{array}
  \end{equation*}
 and  for $s \geq 1$,
  \begin{equation*}
    \begin{array}{rl}
      \mathcal{V} = & \hspace{-4pt} \Big \{ v\ \vert \ v \in L^{\infty}(Q), \ \nabla v \in L^{\infty}(0,T;L^2(\Omega))  \cap L^4(Q) \cap  L^2(0,T;H^1(\Omega)), \\
      & \qquad \Delta v \in L^2(Q), \ \partial_t v \in L^{5/3}(Q) \Big \}.
    \end{array}
  \end{equation*}
  
  We will also need the Banach spaces related to the strong solutions setting:
  \begin{equation*}
    X_p = \{ v \in L^p(0,T;W^{2,p}(\Omega)) \ : \ \partial_t v \in L^p(Q) \}.
  \end{equation*}

\

  \begin{definition}{\bf (Weak solution of \eqref{problema_P_controlado_2})} \label{defi_weak_solution}
    Let $s \geq 1$, $q > 5/2$. Let $f \in L^q(Q)$ and $(u^0, v^0) \in L^p(\Omega) \times W^{2-2/q,q}(\Omega)$, with $p = 1 + \varepsilon$, for some $\varepsilon > 0$, if $s = 1$, and $p = s$, if $s > 1$, be nonnegative functions. A pair $(u,v)$ is called a weak solution of \eqref{problema_P_controlado} if:
    \begin{enumerate}
      \item $u(t,x),v(t,x) \geq 0$ $a.e.$ $(t,x) \in Q$;
      \item $(u,v) \in \mathcal{U} \times \mathcal{V}$;
      \item $u$ and $v$ satisfy the initial conditions (see Remark \ref{remark_initial_conditions_weak_solutions});
      \item the $u$-equation of \eqref{problema_P_controlado} and the boundary condition for $u$ hold in the variational sense
      \begin{equation} \label{eq_u_sol_fraca}
        \int_0^T \int_{\Omega} \partial_t u \ \varphi \ dx \ dt + \int_0^T \int_{\Omega} \nabla u \cdot \nabla \varphi \ dx \ dt = \int_0^T \int_{\Omega} u \nabla v \cdot \nabla \varphi \ dx \ dt,
      \end{equation}
      for all $\varphi \in L^{5s/(4s-3)}(0,T;W^{1,5s/(4s-3)}(\Omega))$, if $s \in [1,2)$, and $\varphi \in L^2(0,T;H^1(\Omega))$, if $s \geq 2$;
      \item  the $v$-equation holds $a.e.$ $(t,x) \in Q$ (in fact,  the $v$-equation is satisfied in $L^{5/3}(Q)$) and, since $\Delta v \in L^2(Q)$, the boundary condition of $v$ in the sense of $H^{-1/2}(\Gamma)$.
    \end{enumerate}
    \hfill $\square$
  \end{definition}

\
  
  \begin{definition}{\bf (Strong solution of \eqref{problema_P_controlado_2})} \label{defi_strong_solution}
    Let $s \geq 1$, $q > 5/2$. Let $f \in L^q(Q)$ and $u^0, v^0 \in W^{2-2/q,q}(\Omega)$ nonnegative functions. A pair $(u,v)$ is called a strong solution of \eqref{problema_P_controlado} if:
    \begin{enumerate}
      \item $u(t,x),v(t,x) \geq 0$ $a.e.$ $(t,x) \in Q$;
      \item $(u,v) \in X_q \times X_q$;
      \item $u$ and $v$ satisfy the initial and boundary conditions of \eqref{problema_P_controlado} (see Remark \ref{remark_initial_and_boundary_conditions_strong_solutions} below);
      \item the $u$-equation and the $v$-equation of \eqref{problema_P_controlado} hold $a.e.$ $(t,x) \in Q$.
    \end{enumerate} 
    \hfill $\square$
  \end{definition}

\

  \begin{remark} \label{remark_initial_and_boundary_conditions_strong_solutions}
    Since $\Delta u, \Delta v \in L^q(Q)$, $u$ and $v$ satisfy the boundary conditions in the sense of $W^{1-1/q,q}(\Gamma)$ (see \cite[Theorem 1.6]{girault2012finite}). Regarding the sense in which the initial condition is attained, we point out that the space $X_p$ is continuously embedded in $C([0, T];W^{2-2/p,p}(\Omega))$ (see \cite[Theorem III.4.10.2]{amann1995linear}).
    \hfill $\square$
  \end{remark}

%
%


\subsection{Optimal control related to weak solutions}
\label{subsection:control weak}

  In this part, the initial conditions of Definition \ref{defi_weak_solution} are supposed to satisfy
  \begin{equation} \label{cond_ini_propriedades}
    (u^0, v^0) \in L^p(\Omega) \times W^{2-2/q,q}(\Omega),
  \end{equation}
  with $p = 1 + \varepsilon$, for some $\varepsilon > 0$, if $s = 1$, and $p = s$, if $s > 1$.
  
  The proof of existence of weak solution to the controlled problem \eqref{problema_P_controlado} is based in the treatment of the uncontrolled model \eqref{problema_P} given in \cite{ViannaGuillen2023uniform}, which is extended  in \cite{guillen2022optimal} to the model \eqref{problema_P_controlado}, to handle the non-smooth control term $f v 1_{\Omega_c}$. Similarly to what was shown in sections \ref{sec:Analytical} and \ref{sec:TimeDiscrete}, an important step in \cite{guillen2022optimal} is obtaining an energy inequality using the change of variable from $(u,v)$ to $(u,z)$, with $z = \sqrt{v + \alpha^2}$, where $\alpha > 0$ is a sufficiently small but fixed real number, independent of $(u,v)$ and $f$. Here, the energy inequality satisfied by the constructed weak solutions of \eqref{problema_P_controlado} will also be written in terms of $(u,z)$. In fact, we consider the energy
  \begin{equation*}
    E(u,z)(t) = \frac{s}{4} \D{\int_{\Omega}}{g(u(t,x)) \ dx} + \frac{1}{2} \int_{\Omega}{\norm{\nabla z(t,x)}{}^2} \ dx,
  \end{equation*}
  where
  \begin{equation*}
    g(u) = \left \{ 
    \begin{array}{rl}
      (u+1)ln(u+1) - u, & \mbox{if } s = 1,  \\[3pt]
      \dfrac{1}{s(s-1)} \, u^s, & \mbox{if } s > 1. 
    \end{array}
    \right.
  \end{equation*}
  The following result of existence of weak solutions to \eqref{problema_P_controlado} is proved in  \cite{guillen2022optimal}.
  \begin{theorem}[\bf Existence of energy inequality weak solutions] \label{teo_existencia_solucao_fraca}
    Given $f \in L^q(Q)$ ($q > 5/2$), there is a non-negative weak solution $(u,v)$ of \eqref{problema_P_controlado} satisfying the following energy inequality
    \begin{equation} \label{desigualdade_integral_limite}
      \begin{array}{l}
        E(u,z)(t_2) + \beta \D{\int_{t_1}^{t_2} \int_{\Omega}}{\norm{\nabla [u + \delta(s)]^{s/2}}{}^2 \ dx} \ dt + \dfrac{1}{4} \D{\int_{t_1}^{t_2} \int_{\Omega}}{u^s \norm{\nabla z}{}^2 \ dx} \ dt \\[6pt]
        + \beta \Big ( \D{\int_{t_1}^{t_2} \int_{\Omega}}{\norm{D^2 z}{}^2 \ dx} \ dt + \D{\int_{t_1}^{t_2} \int_{\Omega}}{\frac{\norm{\nabla z}{}^4}{z^2} \ dx} \ dt \Big ) \\[10pt]
        \leq E(u,z)(t_1) + \mathcal{K}(\norma{f}{L^q(Q)}, \norma{v_0}{W^{2-2/q,q}(\Omega)}),
      \end{array}
    \end{equation}
    for $a.e.$ $t_1,t_2 \in [0,T]$, with $\delta(s) = 1$ if $s \in [1,2)$ and $\delta(s) = 0$ if $s \geq 2$. Here, $\mathcal{K}(\norma{f}{L^q(Q)},\norma{v_0}{W^{2-2/q,q}(\Omega)})$ is a continuous and increasing function with respect to $\norma{f}{L^q(Q)}$ and $\beta > 0$ is a constant, independent of $(u,v,f)$. Moreover, inequality \eqref{desigualdade_integral_limite} is also valid for $t_1 = 0$, with $E(u,z)(0) = E(u^0,v^0)$.    
    Finally, in the case $s > 1$, inequality \eqref{desigualdade_integral_limite} also holds for all $t_2 \in (t_1,T]$.
  \end{theorem}
  Next we introduce the minimization problem. Consider the functional 
  $$J: L^{5s/3}(Q) \times L^2(Q) \times L^q(Q) \longrightarrow \mathbb{R}$$ given by
  \begin{align*}
    & J(u,v,f) : = \dfrac{3\gamma_u}{5s} \int_0^T{\norma{u(t) - u_d(t)}{L^{5s/3}(\Omega)}^{5s/3} \ dt} \\
    & \D + \dfrac{\gamma_v}{2} \int_0^T{\norma{v(t) - v_d(t)}{L^2(\Omega)}^2 \ dt} + \dfrac{\gamma_f}{q} \int_0^T{\norma{f(t)}{L^q(\Omega)}^q \ dt},
  \end{align*}
  where $(u_d,v_d) \in L^{5s/3}(Q) \times L^2(Q)$ represents the desired states, $\gamma_u, \gamma_v > 0$ measure the errors in the states and $ \gamma_f>0$ the cost of the control. In view of the existence result, Theorem \ref{teo_existencia_solucao_fraca}, one could expect the following admissible sets
  \begin{equation*}
    \begin{array}{rl}
      S_{ad}^w = & \hspace{-4pt} \Big \{ (u,v,f) \in X_u \times X_v \times L^q(Q) \ \vert \ (u,v) \mbox{ is a} \\
      & \quad \mbox{weak solution of \eqref{problema_P_controlado} with control } f \}
    \end{array}
  \end{equation*}
  or
  \begin{equation*}
    \begin{array}{rl}
      S_{ad}^E = & \hspace{-4pt} \Big \{ (u,v,f) \in X_u \times X_v \times L^q(Q) \ \vert \ (u,v) \mbox{ is a weak solution of \eqref{problema_P_controlado}} \\
      & \quad \mbox{with control } f \mbox{ and satisfies the energy inequality } \eqref{desigualdade_integral_limite} \Big \}
    \end{array}
  \end{equation*}
  and then state the corresponding minimization problems
  \begin{equation} \label{problema_de_minimizacao_S_ad^w}
    \left \{
    \begin{array}{l} \D
      min \ J(u,v,f) \\
      \mbox{subject to } (u,v,f) \in S_{ad}^w,
    \end{array}
    \right.
  \end{equation}
  or
  \begin{equation} \label{problema_de_minimizacao_S_ad^E}
    \left \{
    \begin{array}{l} \D
      min \ J(u,v,f) \\
      \mbox{subject to } (u,v,f) \in S_{ad}^E.
    \end{array}
    \right.
  \end{equation}
  Thanks to Theorem \ref{teo_existencia_solucao_fraca} we have that both $S_{ad}^w$ and $S_{ad}^E$ are nonempty sets. However, as it is analyzed in \cite{guillen2022optimal} in more detail, we are not able to prove that problem \eqref{problema_de_minimizacao_S_ad^w} or \eqref{problema_de_minimizacao_S_ad^E} has a solution.
  
  Therefore, in order to find an optimal control related to weak solutions of \eqref{problema_P_controlado}, we introduce the following admissible set, for each constant $M>0$:
  \begin{equation*}
    \begin{array}{rl}
      S_{ad}^M = & \hspace{-4pt} \Big \{ (u,v,f) \in X_u \times X_v \times B_q(M) \ \vert \ (u,v) \mbox{ is a weak solution of}  \\
      & \mbox{\eqref{problema_P_controlado} with control } f \mbox{ and satisfies \eqref{desigualdade_integral_limite} changing the} \\
      & \mbox{constant } \mathcal{K}(\norma{f}{L^q(Q)}, \norma{v_0}{W^{2-2/q,q}(\Omega)}) \mbox{ by } \mathcal{K}(M, \norma{v_0}{W^{2-2/q,q}(\Omega)}) \Big \}
    \end{array}
  \end{equation*}
  and the corresponding minimization problem
  \begin{equation} \label{problema_de_minimizacao}
    \left \{
    \begin{array}{l} \D
      min \ J(u,v,f) \\
      \mbox{subject to } (u,v,f) \in S_{ad}^M.
    \end{array}
    \right.
  \end{equation}
  Again, from Theorem \ref{teo_existencia_solucao_fraca}, one has $S_{ad}^M \neq \emptyset$. But now, 
  the following result is proved  in \cite{guillen2022optimal}.
  \begin{theorem}[\bf Existence of optimal control] \label{teo_existencia_controle_otimo}
    For each $M > 0$, the optimal control problem \eqref{problema_de_minimizacao} has at least one global optimal solution, that is, there is $(\overline{u},\overline{v},\overline{f}) \in S_{ad}^M$ such that
    \begin{equation*}
      J(\overline{u},\overline{v},\overline{f}) = \min_{(u,v,f) \in S_{ad}^M} J(u,v,f).
    \end{equation*}
  \end{theorem}
  \begin{remark} \label{remark_energy_structure}
    The existence of weak solutions satisfying an energy inequality is commonly used, for instance, in fluid models,  to prove either weak-strong uniqueness \cite{lions1996incompressible} or large  time behavior \cite{miyakawa1988energy}. In \cite{guillen2022optimal} it is used to prove the existence of global optimal solution. Indeed, considering the minimizing sequence argument, the energy inequality \eqref{desigualdade_integral_limite} is the key point to guarantee that all the possible limits of the minimizing sequence are weak solutions of the controlled model. In fact, the corresponding energy estimates must be strong enough to guarantee that the possible limits of the minimizing sequence are weak solutions of the controlled problem \eqref{problema_P_controlado}. Therefore, if the model does not admit an energy structure, as it seems to be the case in \cite{winkler2015large}, for example, it is not clear how to prove the existence of optimal solution.
    \hfill $\square$
  \end{remark}

\
  
  By construction, one has the following relation between problems \eqref{problema_de_minimizacao_S_ad^w} and \eqref{problema_de_minimizacao}:
  \begin{equation}\label{w-M}
    J_{inf}^w:=\inf_{(u,v,f) \in S_{ad}^w} J(u,v,f) \leq \min_{(u,v,f) \in S_{ad}^M} J(u,v,f).
  \end{equation}
  On the other hand, in \cite{guillen2022optimal}  the following relation between the minimization problems \eqref{problema_de_minimizacao_S_ad^E} and \eqref{problema_de_minimizacao} is obtained for $M$ large enough:
  \begin{theorem} \label{teo_relacao_entre_problemas_de_minimizacao}
    If
    \begin{equation*}
      M^q \geq \frac{q}{\gamma_f} \inf_{(u,v,f) \in S_{ad}^E}{J(u,v,f)},
    \end{equation*}
    one has the inequality
    \begin{equation}\label{M-E}
      \min_{(u,v,f) \in S_{ad}^M} J(u,v,f) \leq \inf_{(u,v,f) \in S_{ad}^E}{J(u,v,f)}:=J_{inf}^E.
    \end{equation}
  \end{theorem}
  To prove the existence of weak solutions satisfying the energy inequality \eqref{desigualdade_integral_limite} for the controlled problem (Theorem \ref{teo_existencia_solucao_fraca}), some procedures used in \cite{ViannaGuillen2023uniform}, and highlighted in Section \ref{sec:Analytical}, are extended in \cite{guillen2022optimal} to the controlled model. First, considering a mollifier regularization of the control $f \in L^q(Q)$ defined via convolution, which results in a sequence  
  (see \cite{brezis2011functional})
  \begin{equation} \label{properties_f_m}
    \begin{array}{c}
      f_m \in C^{\infty}_c(\overline{Q}), \quad \norma{f_m}{L^r(Q)} \leq \norma{f}{L^r(Q)}, \, \mbox{for } r \in [1,q], \\[6pt]
      f_m \rightarrow f \mbox{ strongly in } L^q(Q),
    \end{array}
  \end{equation}
  the truncated controlled problems are defined by
  \begin{equation}
    \left\{
    \begin{array}{l}
      \partial_t u_m - \Delta u_m  = - \nabla \cdot (T^m(u_m) \nabla v_m), \\
      \partial_t v_m - \Delta v_m  = - T^m(u_m)^s v_m + f_m v_m 1_{\Omega_c}, \\
      \partial_\eta u_m \vert_{\Gamma}  = \partial_\eta v_m \vert_{\Gamma} = 0, \quad
      u_m(0)  = u^0_m, \quad v_m(0) = v^0,
    \end{array}
    \right.
    \label{problema_P_m_controlado}
  \end{equation}
  for each $m \in \mathbb{N}$, where the truncation function $T^m \in C^2(\mathbb{R})$ is defined in \eqref{truncamento_limitado_da_identidade}, with $(u^0,v^0)$ satisfying \eqref{cond_ini_propriedades} and $u^0_m \in C^{\infty}(\overline{\Omega})$ being mollifier regularizations of $u^0$ extended to $\mathbb{R}^N$ and having the following properties (see \cite{guillen2022optimal} for more details):
  \begin{equation} \label{convergencia_dado_inicial_u_m}
    u^0_m\geq 0,\quad \int_\Omega u^0_m = \int_\Omega u^0,\quad 
    u^0_m \rightarrow u^0 \mbox{ strongly in } L^p(\Omega).
  \end{equation}
  Then the idea is to pass to the limit in \eqref{problema_P_m_controlado} as $m \to \infty$ to prove the existence of weak solution to \eqref{problema_P_controlado_2} and finish the proof of Theorem \ref{teo_existencia_solucao_fraca} by showing that the weak solution that was found also satisfies the energy inequality \eqref{desigualdade_integral_limite}. Considering first the part of the existence of solution, an important step is to bound the norms $\norma{v_m}{L^{\infty}(Q)}$, uniformly with respect to $m$, and then obtain other $m$-independent estimates that allow us to pass to the limit as $m \to \infty$. The following lemma is essential and is the reason why it is assumed that $f \in L^q(Q)$, with $q > 5/2$.
  \begin{lemma} \label{lema_problema_a_comparar}
    Let $\Omega$ be a bounded domain of $\mathbb{R}^3$ such that $\Gamma$ is of class $C^2$. Let $w^0 \in W^{2-2/q,q}(\Omega)$ and $\tilde{f} \in L^q(Q)$, for some $q > 5/2$. Then the problem
    \begin{equation} \label{problema_a_comparar}
      \left\{\begin{array}{l}
        \partial_t w - \Delta w = \tilde{f} \, w, \quad \mbox{ in } Q, \\
        \partial_\eta w \vert_{\Gamma}  = 0, \quad \mbox{on } (0,T) \times \Gamma, \\
        w(0,x)  = w^0 \  \mbox{in } \Omega,
      \end{array}\right.
    \end{equation}
    has a unique solution
    \begin{equation*}
      w \in C([0,T]; W^{2-2/q,q}(\Omega)) \cap L^q(0,T;W^{2,q}(\Omega)), \ \partial_t w \in L^q(Q),
    \end{equation*}
    In particular, there is a positive constant $C(\|\tilde{f}\|_{L^q(Q)}, \|w^0\|_{W^{2-2/q,q}(\Omega)})$ such that
    \begin{equation} \label{estimate_w_L_infty}
      \norma{w}{L^{\infty}(Q)} \leq C(\|\tilde{f}\|_{L^q(Q)}, \|w^0\|_{W^{2-2/q,q}(\Omega)}).
    \end{equation}
  \end{lemma}

  Similarly to what has been shown in Section \ref{sec:Analytical} ($f \equiv 0$), another important step is to consider the change of variable $z_m(t,x) = \sqrt{v_m(t,x) + \alpha^2}$ and the $(u_m,z_m)$ problem 
    {\small \begin{equation} \label{problema_P_u_m_z_m_controlado}
      \begin{array}{l}
        \partial_t u_m - \Delta u_m = - \nabla \cdot (T^{m}(u_m) \nabla (z_m)^2) \\[6pt]
        \partial_t z_m - \Delta z_m - \dfrac{\norm{\nabla z_m}{}^2}{z_m} = - \dfrac{1}{2} T^{m}(u_m)^s \left(z_m - \dfrac{\alpha^2}{z_m} \right) + \dfrac{1}{2} f_m \left(z_m - \dfrac{\alpha^2}{z_m} \right) 1_{\Omega_c} \\[6pt]
        \partial_\eta u_m  \vert_{\Gamma} =  \partial_\eta z_m  \vert_{\Gamma} = 0 \\[6pt]
        u_m(0) = u^0_m, \quad z_m(0) = \sqrt{v^0 + \alpha^2},
      \end{array}
    \end{equation}}
  which is equivalent to the controlled truncated problem \eqref{problema_P_m_controlado}. Then, using Lemma \ref{lema_problema_a_comparar} the two following equivalent results are proved.

  \begin{theorem} \label{teo_existencia_P_m_controlado}
      Given $f_m$ and $(u^0_m,v^0)$ satisfying \eqref{properties_f_m}, \eqref{convergencia_dado_inicial_u_m} and $v^0 \in W^{2 - 2/q,q}(\Omega)$, respectively, there is a unique solution $(u_m,v_m)$ of \eqref{problema_P_m_controlado} with
      \begin{equation*} 
        \begin{array}{c}
          u_m \in L^{\infty}(0,T;H^1(\Omega)) \cap L^2(0,T;H^2(\Omega)), \ \partial_t u_m \in L^2(0,T;L^2(\Omega)), \\
          v_m \in L^{\infty}(Q) \cap L^{\infty}(0,T;H^2(\Omega)),  \quad \Delta v_m \hbox{ and } \partial_t v_m \in L^2(0,T;H^1(\Omega)),
        \end{array}
      \end{equation*}
      and satisfying
      \begin{equation} \label{positivity_u_v}
        u_m(t,x), v_m(t,x) \geq 0, \ a.e. \ (t,x) \in Q,
      \end{equation}
      \begin{equation} \label{conservation_cells}
        \int_{\Omega}{u_m(t,x) \ dx} = \int_{\Omega}{u^0_m(x) \ dx}= \int_{\Omega}{u^0(x) \ dx}, \ a.e. \ t \in (0,T).
      \end{equation}
      Moreover, there is a positive, continuous and increasing function of $\norma{f}{L^q(Q)}$, $\mathcal{K}_1(\| f\|_{L^q(Q)}, \| v^0\|_{W^{2-2/q,q}(\Omega)})$, also independent of $m$, such that 
      \begin{equation} \label{first_estimates_v_m}
        \norma{v_m}{L^{\infty}(Q)\cap L^2(0,T; H^1(\Omega))} \leq \mathcal{K}_1(\norma{f}{L^q(Q)}, \norma{v^0}{W^{2-2/q,q}(\Omega)}).
      \end{equation}
    \end{theorem}
    \begin{theorem} \label{coro_existencia_P_u_m_z_m_controlado}
      Given $f_m$ and $(u^0_m,v^0)$ satisfying \eqref{properties_f_m}, \eqref{convergencia_dado_inicial_u_m} and $v^0 \in W^{2 - 2/q,q}(\Omega)$, respectively, there is a unique solution $(u_m,z_m)$ of \eqref{problema_P_u_m_z_m_controlado} with
      \begin{equation*} 
        \begin{array}{c}
          u_m \in L^{\infty}(0,T;H^1(\Omega)) \cap L^2(0,T;H^2(\Omega)), \ \partial_t u_m \in L^2(0,T;L^2(\Omega))), \\
          z_m \in L^{\infty}(Q) \cap L^{\infty}(0,T;H^2(\Omega)), \quad \Delta z_m \hbox{ and } \partial_t z_m \in L^2(0,T;H^1(\Omega)),
        \end{array}
      \end{equation*}
      and satisfying the $m$-uniform estimates
      \begin{equation*}
        u_m(t,x) \geq 0 \ \mbox{ and } \ z_m(t,x) \geq \alpha, \ a.e. \ (t,x) \in Q,
      \end{equation*}
      \begin{equation*}
        \int_{\Omega}{u_m(t,x) \ dx} = \int_{\Omega}{u^0_m(x) \ dx} = \int_{\Omega}{u^0(x) \ dx}, \ a.e. \ t \in (0,T),
      \end{equation*}
      \begin{equation} \label{first_estimates_z_m}
        \norma{z_m}{L^{\infty}(Q)}, \norma{z_m}{L^2(0,T;H^1(\Omega))} \leq \mathcal{K}_1(\norma{f}{L^q(Q)}, \norma{v^0}{W^{2-2/q,q}(\Omega)}).
      \end{equation}
    \end{theorem}

    Once we have the existence for the truncated problems, we must obtain $m$-independent estimates and pass to the limit as $m \to \infty$ to prove the existence of weak solutions. This has been done in \cite{guillen2022optimal} by adapting the ideas shown in Section \ref{sec:Analytical} to the controlled problem. First, the following energy inequality is proved:
    \begin{equation} \label{desigualdade_integral_m}
    \begin{array}{l}
      E_m(u_m,z_m)(t_2) + \beta \D{\int_{t_1}^{t_2} \int_{\Omega}}{\norm{\nabla [T^{m}(u_m) + \delta(s)]^{s/2}}{}^2 \ dx} \ dt \\[6pt]
      + \dfrac{1}{4} \D{\int_{t_1}^{t_2} \int_{\Omega}}{T^m(u_m)^s \norm{\nabla z_m}{}^2 \ dx} \ dt \\[6pt]
      + \beta \Big ( \D{\int_{t_1}^{t_2} \int_{\Omega}}{\norm{D^2 z_m}{}^2 \ dx} \ dt + \D{\int_{t_1}^{t_2} \int_{\Omega}}{\frac{\norm{\nabla z_m}{}^4}{z_m^2} \ dx} \ dt \Big ) \\[6pt]
      \leq E_m(u_m,z_m)(t_1) + C (\mathcal{K}_1^2) \D{\int_{t_1}^{t_2} \norma{\nabla z_m}{L^2(\Omega)}^2} \ dt + \mathcal{K}_1^2 \D{\int_{t_1}^{t_2} \norma{f}{L^2(\Omega)}^2} \ dt,
    \end{array}
    \end{equation}
    where $\delta(s) = 1$ if $s \in [1,2)$ and $\delta(s) = 0$ if $s \geq 2$. From \eqref{desigualdade_integral_m}, the $m$-independent estimates were obtained and then used to prove some weak$*$, weak and strong convergences of $(u_m,v_m)$ towards some $(u,v)$ (passing to a subsequence, if necessary) that were finally applied to pass to the limit in \eqref{problema_P_m_controlado}, showing that $(u,v)$ is a weak solution. Additionally, the ideas of \cite{miyakawa1988energy} have been adapted to prove the energy inequality \eqref{desigualdade_integral_limite}. Let
    \begin{equation*}
      E(u,z)(t) = \frac{s}{4} \D{\int_{\Omega}}{g(u(t,x)) \, dx} + \frac{1}{2} \int_{\Omega}{\norm{\nabla z(t,x)}{}^2} dx,
    \end{equation*}
    where
    \begin{equation*}
      g(u) = \left \{ 
      \begin{array}{rl}
      (u+1)ln(u+1) - u, & \mbox{if } s = 1,  \\
      \dfrac{u^s}{s(s-1)}, & \mbox{if } s \in (1,2). 
      \end{array}
      \right.
    \end{equation*}
    It has been proved that $E_m(u_m,z_m) \longrightarrow E(u,z)$ in $L^1(0,T)$, implying
    \begin{equation} \label{convergencia_pontual_E_m}
      E_m(u_m,z_m)(t) \longrightarrow E(u,z)(t) \quad a.e. \ t \in [0,T].
    \end{equation}
    
    With \eqref{convergencia_pontual_E_m} the aforementioned weak$*$, weak and strong convergences were pass to the limit as $m \to \infty$ in the energy inequality \eqref{desigualdade_integral_m}, obtaining \eqref{desigualdade_integral_limite} and concluding the proof of Theorem \ref{teo_existencia_solucao_fraca}.

    Next we give an idea of the proof of Theorems \ref{teo_existencia_controle_otimo} and \ref{teo_relacao_entre_problemas_de_minimizacao}. To prove Theorems \ref{teo_existencia_controle_otimo} a minimizing sequence argument is used. As it was mentioned, $S_{ad}^M$ is non-empty, thanks to Theorem \eqref{teo_existencia_solucao_fraca}. Then there is a sequence $\{ (u_n,v_n,f_n) \} \subset S_{ad}^M$ such that
    \begin{equation*} 
      \lim_{n \to \infty}{J(u_n,v_n,f_n)} = J_{inf} : = \D{\inf_{(u,v,f) \in S_{ad}^M}} J(u,v,f) \geq 0.
    \end{equation*}
    In addition, since $(u_n,v_n,f_n) \in S_{ad}^M$, one has 
  \begin{equation} \label{problema_P_controlado_n}
    \left\{\begin{array}{l}
      \langle \partial_t u_n, \varphi \rangle + \D{\int_\Omega} \nabla u_n\cdot\nabla \varphi \ dx = \D{\int_\Omega} u_n \nabla v_n \cdot \nabla \varphi \ dx \\[6pt]
      \partial_t v_n - \Delta v_n  = - u_n^s v_n + f_n v_n 1_{\Omega_c}, \\[6pt]
      \partial_\eta u_n  \vert_{\Gamma}  =  \partial_\eta v_n  \vert_{\Gamma} = 0, \
      u_n(0)  = u^0, \ v_n(0) = v^0,
    \end{array}\right.
  \end{equation}
   for every $\varphi \in L^{5s/(4s-3)}(0,T;W^{1,5s/(4s-3)}(\Omega))$. Denoting $z_n = \sqrt{v_n + \alpha^2}$, one has  
  \begin{equation} \label{desigualdade_integral_sequencia_minimizante}
    \begin{array}{l}
      E(u_n,z_n)(t_2) + \beta \D{\int_{t_1}^{t_2} \int_{\Omega}}{\norm{\nabla [u_n + \delta(s)]^{s/2}}{}^2  dx} \, dt \\[6pt]
      + \dfrac{1}{4} \D{\int_{t_1}^{t_2} \int_{\Omega}}{u_n^s \norm{\nabla z_n}{}^2  dx} \,dt + \beta \Big ( \D{\int_{t_1}^{t_2} \int_{\Omega}}{\norm{D^2 z_n}{}^2  dx} \, dt \\[6pt]
      + \D{\int_{t_1}^{t_2} \int_{\Omega}}{\frac{\norm{\nabla z_n}{}^4}{z^2} dx} \, dt \Big ) \leq E(u_n,z_n)(t_1) + \mathcal{K}(M,\norma{v_0}{W^{2-2/q,q}(\Omega)}).
    \end{array}
  \end{equation}
  Moreover,
  \begin{equation} \label{limitacao_f_n}
    \norma{f_n}{L^q(Q)} \leq M.
  \end{equation}

  Analogously to the proof of existence of solution, the energy inequality \eqref{desigualdade_integral_sequencia_minimizante} provides $n$-independent bounds that allow us to pass to the limit as $n \to \infty$, proving that there is $(\overline{u},\overline{v},\overline{f}) \in S_{ad}^M$, defined as the limit of a subsequence of $\{ (u_n,v_n,f_n) \}_n$. On the one hand, accounting for the weak lower semi-continuity of the functional $J$, it is possible to prove that
  \begin{equation*}
    J(\overline{u},\overline{v},\overline{f}) \leq \liminf_{n \to \infty}{J(u^n,v^n,f^n)} = J_{inf}.
  \end{equation*}
  On the other hand, since $(\overline{u},\overline{v},\overline{f}) \in S_{ad}^M$, we have $J_{inf} \leq J(\overline{u},\overline{v},\overline{f})$, that is, there is at least one $(\overline{u},\overline{v},\overline{f}) \in S_{ad}^M$ such that $J(\overline{u},\overline{v},\overline{f}) = J_{inf}$.

  To prove Theorem \ref{teo_relacao_entre_problemas_de_minimizacao}, a minimizing sequence argument is also used. Since $S_{ad}^E$ is non-empty, there is a sequence $\{ (u_n,v_n,f_n) \} \subset S_{ad}^E$ such that
  \begin{equation*} 
    \lim_{n \to \infty}{J(u_n,v_n,f_n)} = J_{inf}^E : = \D{\inf_{(u,v,f) \in S_{ad}^E}} J(u,v,f) \geq 0.
  \end{equation*}
  The difference now is that we are not able to prove that the possible limits of $\{ (u_n,v_n,f_n) \}$ are also elements of $S_{ad}^E$. Instead, we are able to prove only that they are elements of $S_{ad}^{M_0}$, for a large enough fixed $M_0 > 0$, yielding the desired result.


\subsection{Optimal control problem subject to strong solutions}
\label{subsection:control strong}

  Now, the initial conditions of the strong solution setting (see Definition \ref{defi_strong_solution}) are supposed to satisfy
  \begin{equation*}
    u^0, v^0 \in W^{2-2/q,q}(\Omega).
  \end{equation*} 
  Moreover, it was necessary to prove in \cite{guillen2023optimal} a regularity criterion that allows one to get existence and uniqueness of global-in-time strong solutions. In fact, the following result has been proved as a first step.
  \begin{theorem} \label{teo_criterio_de_regularidade}
    Let $(u,v)$ be a weak solution of \eqref{problema_P_controlado} with $f \in L^q(Q)$, $q > 5/2$. If, additionally, $u^0,v^0 \in W^{2 - 2/q,q}$ and 
    $$u^s \in L^q(Q),
    $$
     then $v \in X_q$ and $u \in L^{\infty}(Q)$. Moreover $\nabla v \in L^{5q/(5-q)}(Q) \hookrightarrow L^5(Q)$, $u \in X_q$ and $(u,v)$ is the unique strong solution of \eqref{problema_P_controlado}. Finally, there exists $C = C(\norma{u^s}{L^q(Q)},\norma{f}{L^q(Q)},\norma{\nabla u}{L^{5/4}(Q)}) > 0$, which is continuous and increasing with respect to each entry, $\norma{u^s}{L^q(Q)}$, $\norma{f}{L^q}$ and $\norma{\nabla u}{L^{5/4}(Q)}$, such that
    \begin{equation} \label{estimativa_solucao_forte_em_relacao_a_f_nabla_u}
      \norma{(u,v)}{X_q \times X_q} \leq C(\norma{u^s}{L^q(Q)},\norma{f}{L^q(Q)},\norma{\nabla u}{L^{5/4}(Q)}).
    \end{equation}
  \end{theorem}
  \begin{proof}[\bf Idea of the proof.]
    Apply a bootstrapping procedure that yields the desired regularity in a large enough, but finite, number of iterations, which implies inequality \eqref{estimativa_solucao_forte_em_relacao_a_f_nabla_u}. Indeed, a parabolic  regularity result \cite[Theorem 10.22]{feireisl2009singular} is used to prove that $v \in X_q$. Then, a bootstrapping procedure is applied to the $u$-equation of \eqref{problema_P_controlado} to gain regularity for $u$. By hypothesis, we begin with the regularity $u \in L^{p_0}(Q)$, with $p_0 = sq$. Since $q > 5/2$, we have $p_0 > 5/2$ and there is $\alpha > 1$ such that $q = 5 \alpha/2$. Hence, applying Lemma 3.1, we conclude that $u \in X_{qp_0/(q+p_0)}$ and, applying Lemma 2.5 of \cite{guillen2023optimal}, we prove that $u \in L^{p_1}(Q)$, with $p_1 = \alpha p_0$. Repeating this process we increase the regularity of $u$ to $u \in L^{p_n}(Q)$, with $p_n = \alpha^n p_0$, and then Lemma 3.1 of \cite{guillen2023optimal} yields $u \in X_{qp_n/(q+p_n)}$. We do it until we find an index $n_0$ such that $qp_{n_0-1}/(q+p_{n_0-1}) \leq 5/2$ but $qp_{n_0}/(q+p_{n_0}) > 5/2$. In this case $X_{qp_{n_0}/(q+p_{n_0})} \hookrightarrow L^{\infty}(Q)$ and therefore $u \in L^{\infty}(Q)$. Then we finish the proof by applying Lemma 3.1 of \cite{guillen2023optimal} once more.
  \end{proof}
  After, the following result was used to eliminate the dependence of inequality \eqref{estimativa_solucao_forte_em_relacao_a_f_nabla_u} on $\norma{\nabla u}{L^{5/4}(Q)}$.
  \begin{lemma} \label{lema_desigualdade_energia_solucao_forte}
    Assume the existence of the strong solution $(u,v)$ of \eqref{problema_P_controlado} associated to $f \in L^q(Q)$. Let $z = \sqrt{v + \alpha^2}$, for some $\alpha > 0$. Then,
    \begin{equation} \label{limitacao_z_L_infty}
      0 < \alpha \leq z(t,x) \leq \mathcal{K}(\norma{f}{L^q(Q)},\norma{v^0}{W^{2-2/q,q}}).
    \end{equation}
    Moreover, there is $\alpha_0 > 0$, independent of $(u,v,f)$, such that if $0 < \alpha \leq \alpha_0$ then $(u,z)$ satisfies the energy inequality \eqref{desigualdade_integral_limite}.
  \end{lemma}
  Using Lemma \ref{lema_desigualdade_energia_solucao_forte}, 
  it has been proved  in \cite{guillen2023optimal}  that 
  the norm $\norma{\nabla u}{L^{5/4}(Q)}$ can be estimated in terms of $u^0$, $v^0$ and $\norma{f}{L^q(Q)}$. Since the initial data are fixed, the dependence on $(u^0,v^0)$ has been omitted and the following regularity criterion is finally proved.
  \begin{theorem}[\bf Regularity criterion] \label{teo_dependencia_em_relacao_ao_controle}
    Let $(u,v)$ be a weak solution of problem \eqref{problema_P_controlado} with
    $f \in L^q(Q)$. If, additionally, we suppose that $u^s \in L^q(Q)$, then $(u,v) \in X_q \times X_q$ is the unique strong solution of problem \eqref{problema_P_controlado}. Moreover, there is $\hat{\mathcal{K}} = \hat{\mathcal{K}}(\norma{u^s}{L^q(Q)},\norma{f}{L^q(Q)}) > 0$, where $\hat{\mathcal{K}}(\cdot,\cdot)$ is a  continuous and increasing function with respect to each entry, $\norma{u^s}{L^q(Q)}$ and $\norma{f}{L^q}$, such that
    \begin{equation} \label{estimativa_solucao_forte_em_relacao_a_f}
      \norma{(u,v)}{X_q \times X_q} \leq \hat{\mathcal{K}}(\norma{u^s}{L^q(Q)},\norma{f}{L^q(Q)}).
    \end{equation}
  \end{theorem}
  \begin{remark}
    Differently from the regularity criteria given in \cite{guillen2020optimal} and \cite{lopez2021optimal} for other chemotaxis models, for example, the regularity criterion given in Theorem~\ref{teo_dependencia_em_relacao_ao_controle} is sharp because hypothesis $f \in L^q(Q)$ with $q > 5/2$ also appears in the weak solution setting, being essential in the proof given in \cite[Lemma 3.1]{guillen2022optimal} of existence of weak solutions of \eqref{problema_P_controlado}.
    \hfill $\square$
  \end{remark}

\

  \begin{remark}
    Since $q > 5/2$, we have $X_q \hookrightarrow L^{\infty}(Q)$. Then Theorem \ref{teo_dependencia_em_relacao_ao_controle} also gives a regularity hypothesis over a weak solution of the controlled problem \eqref{problema_P_controlado} avoiding blow up at finite time.
    \hfill $\square$
  \end{remark}

\

  Next, the following optimal control problem can be introduced. Let $\mathcal{F}$ be a closed and convex subset of $L^q(Q)$, for a given $q > 5/2$. Consider the cost functional $J: L^{sq}(Q) \times L^2(Q) \times \mathcal{F} \longrightarrow \mathbb{R}$ given by
  \begin{equation} \label{funcional_J}
    \begin{array}{l}
      J(u,v,f) : = \dfrac{\gamma_u}{s q} \D \int_0^T{\norma{u(t) - u_d(t)}{L^{sq} }^{sq} \ dt} \\[6pt]
      \D + \dfrac{\gamma_v}{2} \int_0^T{\norma{v(t) - v_d(t)}{L^2 }^2 \ dt} + \dfrac{\gamma_f}{q} \int_0^T{\norma{f(t)}{L^q(\Omega_c)}^q \ dt},
    \end{array}
  \end{equation}
  where $(u_d,v_d) \in L^{sq}(Q) \times L^2(Q)$ represents the desired states and the parameters $\gamma_u, \gamma_v, \gamma_f \geq 0$ measure the costs of the states and control. In addition, we assume
  \begin{equation} \label{hipotese_sobre_os_custos_gamma}
    \begin{array}{c}
      \gamma_u > 0 \ \mbox{ and} \\
      \gamma_f > 0 \ \mbox{ or } \ \mathcal{F} \mbox{ is bounded in } L^q(Q).
    \end{array}
  \end{equation}
  
  Note that the functional $J(u,v,f)$ is well defined for weak solutions of the controlled problem with $f,u^s \in L^q(Q)$ but, because of the regularity criterion given in Theorem~\ref{teo_dependencia_em_relacao_ao_controle}, such solutions are actually the strong solutions. This allows us to define the optimal control problem as
  \begin{equation} \label{problema_de_minimizacao_forte1}
    min \ J(u,v,f) \mbox{ subject to } (u,v,f) \in S_{ad},
  \end{equation}
  with
  \begin{equation*}
    \begin{array}{rcl}
      S_{ad}  =  \{ (u,v,f) \in X_q \times X_q \times \mathcal{F} \ \vert \ (u,v) \mbox{ the strong solution of \eqref{problema_P_controlado}} \}.
    \end{array}
  \end{equation*}
  
  Since given $f \in \mathcal{F}$ one cannot assure in general the existence of a strong solution $(u,v)$ associated to $f$, we fix the hypothesis
  \begin{equation*} \label{hipotese_conjunto_admissivel_nao_vazio}
    S_{ad} \neq \emptyset.
  \end{equation*}
  As we mentioned in the introductory text of Section \ref{sec:Control}, it is possible to prove that $S_{ad} \neq \emptyset$ if $\Omega_c = \Omega$ or if $\Omega$ is a two dimensional domain. We have the following.
  
  \begin{theorem}[\bf Existence of optimal control] \label{teo_existencia_controle_otimo_forte}
    Assuming $S_{ad} \neq \emptyset$, the optimal control problem \eqref{problema_de_minimizacao_forte1} has at least one global optimal solution $(\overline{u},\overline{v},\overline{f}) \in S_{ad}$.
  \end{theorem}

  Similarly to Theorem \ref{teo_existencia_controle_otimo}, Theorem \ref{teo_existencia_controle_otimo_forte} was also proved by means of a minimizing sequence argument \cite{guillen2023optimal}. Indeed, assuming $S_{ad} \neq \emptyset$, we conclude that there is a sequence $\{ (u_n,v_n,f_n) \}_n \subset S_{ad}$ such that
  \begin{equation*}
    \lim_{n \to \infty}{J(u_n,v_n,f_n)} = \inf_{(u,v,f) \in S_{ad}}{J(u,v,f)} \geq 0.
  \end{equation*}
  By using the definition of $J$,  $\norma{u_n}{L^{sq}(Q)}$ and $\norma{f_n}{L^q(Q)}$ are bounded and, by inequality \eqref{estimativa_solucao_forte_em_relacao_a_f} of Theorem~\ref{teo_dependencia_em_relacao_ao_controle}, we also conclude that
  \begin{equation} \label{limitacao_seq_minimizante_forte}
    \norma{(u_n,v_n)}{X_q \times X_q} \mbox{ is bounded}.
  \end{equation}
  This bound is sufficient to prove that all the possible limits of $\{ (u_n,v_n,f_n) \}_n$ belong to $S_{ad}$ and to conclude the proof of Theorem \ref{teo_existencia_controle_otimo_forte}.

\

  \begin{remark}
    Note that the passage from inequality \eqref{estimativa_solucao_forte_em_relacao_a_f_nabla_u} to \eqref{estimativa_solucao_forte_em_relacao_a_f}, eliminating the dependence on $\norma{\nabla u}{L^{5/4}(Q)}$, has been very important to get \eqref{limitacao_seq_minimizante_forte}. In fact, if we had only \eqref{estimativa_solucao_forte_em_relacao_a_f_nabla_u} and wanted to prove existence of optimal control, we would need to include the norm $\norma{\nabla u - \nabla u_d}{L^{5/4}(Q)}$ in the definition of the functional $J$. However, it would mean to require the desired state $u_d$ to have first order derivatives, which can be too restrictive from the applications point of view.
    \hfill $\square$
  \end{remark}

\
  
  After proving existence of optimal solution, first order optimality conditions for any local optimum were established in \cite{guillen2023optimal}, via a generic Lagrange multipliers theorem, see Theorem \ref{teo_existencia_multiplicadores_de_Lagrange} below. To apply this generic result, it is necessary to write the optimal control problem \eqref{problema_de_minimizacao_forte1} in an abstract form. Then, we introduce the Banach spaces $\mathbb{X} = \widetilde{X}_q \times \widetilde{X}_q \times L^q(Q),$ and $\mathbb{Y} = L^q(Q) \times L^q(Q)$  where $\widetilde{X}_q = \{ w \in X_q \ \vert \ \partial_{\boldsymbol{n}} w \vert_{\Gamma} = 0 \}$, and the operator $G = (G_1,G_2): \mathbb{X} \rightarrow \mathbb{Y}$ defined for each $r = (u,v,f) \in \mathbb{X}$ as
  \begin{equation*}
    \left \{
    \begin{array}{rl}
      G_1(r) & = \partial_t u - \Delta u + \nabla \cdot (u \nabla v) \\
      G_2(r) & = \partial_t v - \Delta v + u^sv - fv \ 1_{\Omega_c}.
    \end{array}
    \right.
  \end{equation*}
  By using hypothesis $S_{ad}\not=\emptyset$,  there exists $(\hat{u},\hat{v},\hat{f})\in S_{ad}$. Then, we consider the space
  \begin{equation*}
    \widehat{X}_q = \{ w \in \widetilde{X}_q \ \vert \ w(0,x) = 0 \}
  \end{equation*}
  and we define $\mathbb{M}$, a closed and convex subset of $\mathbb{X}$, as
  $$\mathbb{M} = (\widehat{u}, \widehat{v}, \widehat{f}) + \widehat{X}_q \times \widehat{X}_q \times (\mathcal{F} - \widehat{f}).$$
  Therefore, we rewrite the optimal control problem \eqref{problema_de_minimizacao_forte1} as
  \begin{equation} \label{problema_de_minimizacao_abstrato}
    \min_{r \in \mathbb{M}} J(r) \ \mbox{ subject to } \ G(r) = 0,
  \end{equation}
  where $J: \mathbb{X} \rightarrow \mathbb{R}$, $G: \mathbb{X} \rightarrow \mathbb{Y}$, with $\mathbb{X}$ and $\mathbb{Y}$ being Banach spaces and $\mathbb{M}$ a closed convex subspace of $\mathbb{X}$. Note that we can rewrite the admissible set as $S_{ad} = \{ r \in \mathbb{M} \ \vert \ G(r) = 0 \}$.

  We must also define the concepts of Lagrangian functional, Lagrange multiplier and regular point.
  \begin{definition}{\bf (Lagrangian)}
    The functional $\mathcal{L}: \mathbb{X} \times \mathbb{Y}' \rightarrow \mathbb{R}$, given by
    \begin{equation} \label{Lagrangian_functional}
      \mathcal{L}(r,\xi) = J(r) - \dist{\xi}{G(r)}_{\mathbb{Y}'},
    \end{equation}
    is called the Lagrangian functional related to problem \eqref{problema_de_minimizacao_abstrato}.
    \hfill $\square$
  \end{definition}

\

  \begin{definition}{\bf (Lagrange multipliers)}
    Let $\overline{r} \in S$ be a local optimal solution of problem \eqref{problema_de_minimizacao_abstrato}. Suppose that $J$ and $G$ are Fréchet differentiable in $\overline{r}$, the derivatives being denoted by $J'(\overline{r})$ and $G'(\overline{r})$, respectively. Then, $\xi \in \mathbb{Y}'$ is called a Lagrange multiplier for \eqref{problema_de_minimizacao_abstrato} at the point $\overline{r}$ if
    \begin{equation} \label{expressao_multiplicador_de_Lagrange}
      \mathcal{L}'(\overline{r},\xi)[c] = J'(\overline{r})[c] - \dist{\xi}{G'(\overline{r})[c]}_{\mathbb{Y}'} \geq 0, \ \forall c \in \mathcal{C}(\overline{r}),
    \end{equation}
    where $\mathcal{C}(\overline{r}) = \{ \theta (r - \overline{r}) \ \vert \ r \in \mathbb{M}, \ \theta \geq 0 \}$ is the conical hull of $\overline{r} \in \mathbb{M}$.
    \hfill $\square$
  \end{definition}

\

  \begin{definition}{\bf (Regular point)} \label{defi_regular_point}
    A point $\overline{r} \in \mathbb{M}$ is called a regular point if $G'(\overline{r})[\mathcal{C}(\overline{r})] = \mathbb{Y}$.
    \hfill $\square$
  \end{definition}

\

  Finally, we state the theorem on the existence of Lagrange multipliers.
  \begin{theorem}{\bf (\cite{zowe1979regularity})} \label{teo_existencia_multiplicadores_de_Lagrange}
    Let $\overline{r} \in S$ be a local optimal solution of problem \eqref{problema_de_minimizacao_abstrato}. Suppose that $J$ is Fréchet differentiable and $G$ is continuously Fréchet differentiable. If $\overline{r}$ is a regular point, then there exists a  Lagrange multiplier for problem \eqref{problema_de_minimizacao_abstrato} at $\overline{r}$.
  \end{theorem}

  Accounting for the hypotheses of Theorem \ref{teo_existencia_multiplicadores_de_Lagrange} and \eqref{expressao_multiplicador_de_Lagrange}, we need to prove that $J$ and $G$ are Fréchet differentiable and find the expression of the derivatives. Let
  \begin{equation} \label{rhs-adjoint}
     g_\lambda= \gamma_u sgn(\overline{u} - u_d) \norm{\overline{u} - u_d}{}^{sq - 1}
     \quad  \hbox{and} \quad
     g_\eta=\gamma_v (\overline{v} - v_d).
   \end{equation}
   We have the following results.

  \begin{lemma}[\bf \cite{guillen2023optimal}] \label{lema_J_diferenciavel}
    The functional $J: \mathbb{X} \rightarrow \mathbb{R}$ is Fréchet differentiable and the Fréchet derivative of $J$ in $\overline{r} = (\overline{u},\overline{v},\overline{f}) \in \mathbb{X}$ in the direction $c = (U,V,F) \in \mathbb{X}$ is 
    \begin{equation} \label{derivada_Frechet_J}
        J'(\overline{r})[c] 
        = \D{\int_0^T \int_{\Omega}} (g_\lambda\, U + g_\eta\, V) dx \ dt
         + \gamma_f \D{\int_0^T \int_{\Omega_c}} sgn(\overline{f}) \norm{\overline{f}}{}^{q-1} F \ dx \ dt,
    \end{equation}
    where $g_\lambda,g_\eta$ are defined in \eqref{rhs-adjoint}.
  \end{lemma}
  
  \begin{lemma}[\bf \cite{guillen2023optimal}] \label{lema_G_diferenciavel}
    The operator $G: \mathbb{X} \rightarrow \mathbb{Y}$ is continuously Fréchet differentiable and the Fréchet derivative of $G$ in ${r} = ({u},{v},{f}) \in \mathbb{X}$ in the direction $c = (U,V,F) \in \mathbb{X}$ is the operator $G'({r})[c] = (G_1'({r})[c],G_2'({r})[c])$ given by
    \begin{equation} \label{derivada_Frechet_G}
    \left \{
    \begin{array}{rl}
      G_1'({r})[c] & \hspace{-2mm} = \partial_t U - \Delta U + \nabla \cdot (U \nabla{v}) + \nabla \cdot ({u} \nabla V) \\
      G_2'({r})[c] & \hspace{-2mm} = \partial_t V - \Delta V + s\, {u}^{s-1} U{v} +{u}^s V -{f} V \ 1_{\Omega_c} - F{v} \ 1_{\Omega_c}.
    \end{array}
    \right.
    \end{equation}
  \end{lemma}

\

  Considering \eqref{expressao_multiplicador_de_Lagrange} and the expressions of the derivatives of $J$ and $G$, it has been proved in \cite{guillen2023optimal} that the existence of a Lagrange multiplier $(\lambda,\eta)$ associated to a local optimum $(\overline{u}, \overline{v}, \overline{f})$ is equivalent to finding a pair $(\lambda, \eta)$ which is a very weak solution of the problem 
  \begin{equation} \label{problema_adjunto_ao_linearizado}
     \left \{
     \begin{array}{l}
       - \partial_t \lambda - \Delta \lambda - \nabla \overline{v} \cdot \nabla \lambda + s \overline{u}^{s-1} \overline{v} \eta = g_\lambda, \\[6pt]
       - \partial_t \eta - \Delta \eta + \overline{u}^s \eta - \overline{f} \eta \ 1_{\Omega_c} + \nabla \cdot (\overline{u} \nabla \lambda) = g_\eta, \\[6pt]
       \partial_{\boldsymbol{n}} \lambda \vert_{\Gamma} = \partial_{\boldsymbol{n}} \eta \vert_{\Gamma} = 0, \ \lambda(T,x) = \eta(T,x) = 0,
     \end{array}
     \right.
   \end{equation}
   in the sense of the following definition.

\
  
  \begin{definition}{\bf (Very weak solution of \eqref{problema_adjunto_ao_linearizado})} \label{defi_very_weak_solution}
     Let $s \geq 1$, $q > 5/2$ and $q' = q/(q - 1)$. A pair $(\lambda,\eta) \in L^{q'}(Q) \times L^{q'}(Q)$ is called a very weak solution of \eqref{problema_adjunto_ao_linearizado} if $(\lambda,\eta)$ satisfies \eqref{problema_adjunto_ao_linearizado} in the sense of the dual space of $X_q \times X_q$, that is, the following variational formulation holds for any $U,V \in X_q$ with 
     $\partial_{\boldsymbol{n}} U  \vert_{\Gamma} = \partial_{\boldsymbol{n}} V  \vert_{\Gamma}=0$ and $U(0)=V(0)=0$:
     {\small \begin{equation}\label{adjoint-1}
      \begin{array}{l}
      \D  {\int_0^T \!\!\int_{\Omega}} \lambda \Big ( \partial_t U - \Delta U + \nabla \cdot (U \nabla \overline{v}) \Big ) 
      + \int_0^T \!\! \int_{\Omega} s \overline{u}^{s-1} \overline{v} \eta \ U 
      = \int_0^T \int_{\Omega}g_\lambda U ,
      \end{array}
    \end{equation}}
    {\small \begin{equation}\label{adjoint-2}
      \begin{array}{l}
   \D    {\int_0^T\!\! \int_{\Omega}} \eta \Big ( \partial_t V - \Delta V + \overline{u}^s V - \overline{f} V 1_{\Omega_c} \Big ) 
   + \int_0^T \!\!\int_{\Omega} \lambda \ \nabla \cdot (\overline{u} \nabla V) 
   = \D\int_0^T \!\!\int_{\Omega} g_\eta  V .
     \end{array}
    \end{equation}}
     \hfill $\square$
   \end{definition}
   The last hypothesis to be proved in order to apply Theorem \ref{teo_existencia_multiplicadores_de_Lagrange} is that a local optimal solution is a regular point. In fact, in \cite{guillen2023optimal} it is proved that any $(\overline{u},\overline{v},\overline{f}) \in S_{ad}$ is a regular point, because for any $(g_U,g_V) \in \mathbb{Y}$, there is $ (U,V) \in \widehat{X}_q \times \widehat{X}_q$ such that
  \begin{equation} \label{problema_linearizado}
    \left \{
    \begin{array}{l}
      \partial_t U - \Delta U = - \nabla \cdot (U \nabla \overline{v}) - \nabla \cdot (\overline{u} \nabla V)  + g_U \\
      \partial_t V - \Delta V = - s \overline{u}^{s-1} U \overline{v} - \overline{u}^s V + \overline{f} V \ 1_{\Omega_c} + g_V.
    \end{array}
    \right.
  \end{equation}

\

  \begin{remark}
    Problem \eqref{problema_linearizado} is called the linearized problem related to \eqref{problema_P_controlado}. Note that the Lagrange multiplier system \eqref{problema_adjunto_ao_linearizado} is the adjoint of the linearized problem \eqref{problema_linearizado}.
    \hfill $\square$
  \end{remark}

\

  To prove the  existence of solution to the linearized problem \eqref{problema_linearizado}, the following general prototype of a linearized problem related to chemotaxis models has been studied in \cite{guillen2023optimal}:
  \begin{equation} \label{problema_linearizado_geral}
    \left \{
    \begin{array}{l}
      \partial_t U - \Delta U + a_1 U + b_1 V + \nabla \cdot (U \Vec{c}_1) + \nabla \cdot (d \nabla V) = g_U, \\
      \partial_t V - \Delta V + a_2 V + b_2 U + \Vec{c}_2 \cdot \nabla V = g_V, \\
      \partial_{\boldsymbol{n}} U \vert_{\Gamma} = \partial_{\boldsymbol{n}} V \vert_{\Gamma} = 0, \ U(0,x) = V(0,x) = 0,
    \end{array}
    \right.
  \end{equation}
  where the coefficient functions $a_i, b_i, d$ and $\Vec{c}_i$ are data defined in $Q$. Considering the Banach space
  \begin{equation*}
    W_2 = \{ v \in L^2(H^{1}) \ : \ \partial_t v \in L^{2}((H^1)') \},
  \end{equation*}
  the following general result was proved.
  \begin{theorem}[\bf \cite{guillen2023optimal}] \label{teo_problema_linearizado_geral}
     Let $a_i \in L^{5/2}(Q)$ and $\Vec{c}_i \in L^5(Q)^3$ with $\Vec{c}_i \cdot \Vec{n} \vert_{\Gamma} = 0$ for $i = 1, 2$.
    \begin{enumerate}
      \item if $b_i \in L^{5/2}(Q)$ and $d \in L^{\infty}(Q)$ we have:
      \begin{enumerate}
        \item if $g_U,g_V \in L^{10/7}(Q)$ then there is a weak solution $(U,V) \in W_2 \times W_2$ of \eqref{problema_linearizado_geral};
        \item if $g_U,g_V \in L^{10/9}(Q)$ and $\nabla d \in L^5(Q)^3$ then there is a very weak solution $(U,V) \in L^2(Q) \times L^2(Q)$ of \eqref{problema_linearizado_geral};
      \end{enumerate}
      \item if $b_1 \in L^{5/3}(Q)$ and $b_2, d \in L^5(Q)$ we have:
      \begin{enumerate}
        \item if $g_U \in L^{10/7}(Q)$ and $g_V \in L^2(Q)$ then there is a weak-strong solution $(U,V) \in W_2 \times X_2$ of \eqref{problema_linearizado_geral};
        \item if $g_U \in L^{10/9}(Q)$ and $g_V \in L^{10/7}(Q)$ then there is a very weak-weak solution $(U,V) \in L^2(Q) \times W_2$ of \eqref{problema_linearizado_geral}.
      \end{enumerate}
    \end{enumerate}
  \end{theorem}
   We are now in position to introduce the first order optimality conditions.
   \begin{theorem}[\bf \cite{guillen2023optimal}] \label{teo_existencia_multiplicador_de_Lagrange_aplicado}
    Assume \eqref{assumptions_set} and let $(\overline{u},\overline{v},\overline{f}) \in S_{ad}$ be a local optimal solution of \eqref{problema_de_minimizacao_forte1}. Then there exists a unique Lagrange multiplier $(\lambda,\eta) \in L^{q'}(Q) \times L^{q'}(Q)$ which is a very weak solution of the optimality system \eqref{problema_adjunto_ao_linearizado} and the following optimality condition holds:
    \begin{equation} \label{condicao_otimalidade_fraco_F}
      \D{\int_0^T \int_{\Omega_c}} (\gamma_f sgn(\overline{f}) \norm{\overline{f}}{}^{q-1} + \overline{v}\, \eta)(f - \overline{f}) \ dx \ dt \geq 0, \quad \forall f \in \mathcal{F}.
    \end{equation}
  \end{theorem}

\

  \begin{remark}
    If $\gamma_f > 0$ and there is no convex constraint on the control, that is $\mathcal{F} = L^q(Q)$, then \eqref{condicao_otimalidade_fraco_F} is equivalent to
    $\gamma_f sgn(\overline{f}) \norm{\overline{f}}{}^{q-1} + \overline{v} \,\eta = 0$.
    Since $\overline{v} \geq 0$, we conclude the following explicit expression for the control $\overline{f} = - sgn(\eta) \left ( \dfrac{1}{\gamma_f} \overline{v}\, \norm{\eta}{} \right )^{1/(q-1)}$.
    \hfill $\square$
  \end{remark}

\

  To finish, the regularity of the Lagrange multiplier $(\lambda,\eta)$ is also studied. The general linear system \eqref{problema_linearizado_geral} and Theorem \ref{teo_problema_linearizado_geral} were useful once again, this time to prove the following result.
  \begin{theorem}[\bf \cite{guillen2023optimal}] \label{teo_regularidade_adicional_multiplicadores_de_Lagrange}
    Assume \eqref{assumptions_set} and let $(\overline{u},\overline{v},\overline{f}) \in S_{ad}$ be a local optimal of problem \eqref{problema_de_minimizacao_forte1}. It holds:
    \begin{enumerate}
      \item if $g_{\lambda} \in L^p(Q)$, for $p\in [10/9 , 10/7)$, then 
      the Lagrange multiplier $(\lambda,\eta) \in L^2(Q) \times L^2(Q)$ and satisfies  \eqref{problema_adjunto_ao_linearizado} in the very weak sense (as in \eqref{adjoint-1}-\eqref{adjoint-2});
      \item if $g_{\lambda} \in L^p(Q)$, for $p\in [10/7 , 2]$, then 
      the Lagrange multiplier $(\lambda,\eta) \in X_p \times X_p$ and satisfies 
      \eqref{problema_adjunto_ao_linearizado} in the strong sense, that is, $a.e.$ $(t,x)$ in $Q$.
    \end{enumerate}
  \end{theorem}

  \begin{remark}
    Since $v_d \in L^2(Q)$, which implies $g_{\eta} \in L^2(Q)$, the previous analysis for $p > 2$ does not seem to lead to  more relevant conclusions.
    \hfill $\square$
  \end{remark}

\

  \begin{remark}
    To guarantee that the terms of the functional $J$ given in \eqref{funcional_J} make sense, it is enough that $u_d \in L^{\tilde q}(Q)$, with $\tilde{q} \geq sq$, and $v_d \in L^2(Q)$. With this regularity,  $g_{\eta} \in L^2(Q)$ and $g_{\lambda} \in L^p(Q)$, for a power $p = p(s,q,\tilde{q}) = \tilde{q}/(sq - 1)$. Hence the regularity of $g_{\lambda}$ depends on $s \geq 1$, $q > 5/2$ and $\tilde{q} \geq sq$, and is decreasing with respect to $s$, with $p(s,q,\tilde{q}) \to 1$ as $s \to \infty$. For instance if  $\tilde{q} = sq$, we have $p = sq/(sq - 1)$. In this case, since $s \geq 1$ and $q > 5/2$, then $p \in (1,5/3)$. Let us fix $q > 5/2$ close to $5/2$ and vary the values of $s$. Then,  if $s \in [1,10/3q]$  we are in the case $2$ of Theorem \ref{teo_regularidade_adicional_multiplicadores_de_Lagrange}, and if $s \in (10/3q,10/q]$  we are in the case $1$ of Theorem \ref{teo_regularidade_adicional_multiplicadores_de_Lagrange}. But, if $s > 10/q$ then $p \in (1,10/9)$, hence  Theorem \ref{teo_regularidade_adicional_multiplicadores_de_Lagrange} doesn't give  
 additional regularity for the Lagrange multiplier.
    \hfill $\square$
  \end{remark}


\bmhead{Acknowledgments}

This work was partially funded by grant PGC2018-098308-B-I00 (MCI/AEI/FEDER, UE). The second author has also been financed in part by grant US-1381261 (US/JUNTA/FEDER, UE) and grant P20-01120 (PAIDI/JUNTA/FEDER, UE).

\section*{Declarations}

\subsection*{Competing interests}

    The authors have no relevant financial or non-financial interests to disclose.

  \subsection*{Author contributions}

    All authors contributed to the study conception and design. All authors read and approved the final manuscript.




\nocite{label}
\bibliography{sn-bibliography}


\end{document}